\documentclass[11pt,reqno]{amsart}

\usepackage{amsmath, amsfonts, amssymb, amscd, enumerate}
\theoremstyle{plain}
\newtheorem{theo}{Theorem}[section]
\newtheorem{lem}[theo]{Lemma}
\newtheorem{prop}[theo]{Proposition}

\theoremstyle{definition}
\newtheorem{rem}[theo]{Remark}
\newtheorem{example}[theo]{Example}
\newtheorem{definition}[theo]{Definition}
\newenvironment{pf}{\noindent{\it Proof. }}{$\square$\par\medskip}

\theoremstyle{plain}

\theoremstyle{definition}

\renewcommand{\=}{\overset{\operatorname{def}}{=}}

\newcommand{\beq}{\begin{equation}}
\newcommand{\eeq}{\end{equation}}
\renewcommand{\a}{\alpha}
\renewcommand{\b}{\beta}

\renewcommand{\d}{\delta}

\newcommand{\f}{\varphi}
\newcommand{\g}{\gamma}
\newcommand{\h}{\eta}

\renewcommand{\l}{\lambda}
\renewcommand{\o}{\omega}

\renewcommand{\r}{\rho}

\renewcommand{\t}{\tau}

\newcommand{\z}{\zeta}
\newcommand{\D}{\Delta}

\newcommand{\G}{\Gamma}

\renewcommand{\O}{\Omega}

\newcommand{\bA}{\mathbb{A}}

\newcommand{\bC}{\mathbb{C}}

\newcommand{\bF}{\mathbb{F}}
\newcommand{\bG}{\mathbb{G}}
\newcommand{\bR}{\mathbb{R}}

\newcommand{\bJ}{\mathbb{J}}

\newcommand{\gC}{\mathfrak{C}}
\newcommand{\gD}{\mathfrak{D}}

\newcommand{\cA}{\mathcal{A}}
\newcommand{\cB}{\mathcal{B}}
\newcommand{\cC}{\mathcal{C}}

\newcommand{\cF}{\mathcal{F}}

\newcommand{\cH}{\mathcal{H}}

\newcommand{\cL}{\mathcal{L}}

\newcommand{\cN}{\mathcal{N}}

\newcommand{\cU}{\mathcal{U}}
\newcommand{\cV}{\mathcal{V}}

\newcommand{\cZ}{\mathcal{Z}}

\newcommand{\Jst}{J_{\operatorname{st}}}
\newcommand{\dcst}{d^c_{\operatorname{st}}}

\newcommand{\HH}{\cH ess}

\renewcommand{\square}{\kern1pt\vbox
{\hrule height 0.6pt\hbox{\vrule width 0.6pt\hskip 3pt
\vbox{\vskip 6pt}\hskip 3pt\vrule width 0.6pt}\hrule height0.6pt}\kern1pt}

\DeclareMathOperator{\Span}{Span}

\renewcommand\Re{\operatorname{Re}}
\renewcommand\Im{\operatorname{Im}}
\newcommand\Hom{\operatorname{Hom}}

\newcommand{\wt}{\widetilde}
\newcommand{\wh}{\widehat}

\newcommand{\Psh}{\operatorname{Psh}}

\newcommand{\be}{\begin{equation}}
\newcommand{\ee}{\end{equation}}

\def\<#1,#2>{\langle\,#1,\,#2\,\rangle}
\newcommand{\arr}{\begin{array}{rlll}}
\newcommand{\ea}{\end{array}}
\newcommand{\bea}{\begin{eqnarray}}
\newcommand{\eea}{\end{eqnarray}}
\newcommand{\bean}{\begin{eqnarray*}}
\newcommand{\eean}{\end{eqnarray*}}

\def\sideremark#1{\ifvmode\leavevmode\fi\vadjust{
\vbox to0pt{\hbox to 0pt{\hskip\hsize\hskip1em
\vbox{\hsize3cm\tiny\raggedright\pretolerance10000
\noindent #1\hfill}\hss}\vbox to8pt{\vfil}\vss}}}

\newcounter{ssig}
\newcounter{ttig}
\renewcommand{\exp}{\Phi}

\title[Stationary  disks and Green functions]
{Stationary  disks  and Green functions \\ in almost complex domains}

\author{G. Patrizio and A. Spiro}

\begin{document}

\begin{abstract} 
Using  generalized Riemann maps,  normal forms for almost complex domains $(D, J)$ with singular foliations  by stationary disks are defined. Such  normal forms are used to construct counterexamples and to determine intrinsic conditions,  under which the stationary disks  are  extremal disks  for the Kobayashi metric or determine solutions to  almost complex Monge-Amp\`ere equation.

  \end{abstract}


\thanks{{\it Acknowledgments}. This research was partially supported by the Project MIUR ÒGeometric Properties of Real and Complex ManifoldsÓ,  Project MIUR ÒDifferential Geometry and Global AnalysisÓ and by GNSAGA of INdAM}

\address{
Dipartimento di Matematica ``U. Dini'', 
Universit\`a di Firenze, 
Firenze, 
ITALY}
\email{patrizio@math.unifi.it}
\address{
Scuola di Scienze e Tecnologie, Universit\`a di Camerino, 
 Camerino, 
ITALY}
\email{andrea.spiro@unicam.it}
\maketitle

\null \vspace*{-.25in}

\section{Introduction}
\setcounter{section}{1}
\setcounter{equation}{0}

%
%
%
%
%
%
%

Let $D \subset M$ be a domain of an almost complex manifold $(M, J)$ with smooth boundary $\partial D$ and  $\wt \bJ$ the canonical almost complex structure of  $T^*M$ determined by $J$.  We recall that a smooth, proper $J$-holomorphic embedding $f:  \D \to D$ of  the unit disk  into $D$
 is called {\it stationary disk} if  there exists a map $\wt f: \overline \D \to T^* M\setminus \{\text{zero section}\}$, which is $\wt \bJ$-holomorphic, projects onto $f$ and is such that, for any $\z \in \partial \D$, the 1-form $\z^{-1} \cdot  \wt f(\z) \in T^*_{f(\z)}M $ vanishes identically on $T_{f(\z)}\partial D$ (see \S \ref{preliminaries}). \par
The   stationary disks of almost complex domains have been introduced  by Coupet, Gaussier and Sukhov in \cite{CGS, CGS1} (see also \cite{CGS2}). They are useful biholomorphic invariants  of almost complex  domains and  constitute a natural generalization of the   stationary disks of strictly linearly  convex domains  of $\bC^n$, considered for the first time  in  celebrated Lempert's papers  on extremal disks and Kobayashi metrics  (\cite{Le, Le2}).\par
\smallskip
Existence and uniqueness results on stationary disks,  with prescribed center and direction,  has been established in various contexts, both in the integrable and non-integrable case (see e.g. \cite{Le, Pa, Tu, CGS, SS, PS1, CGS2}).  
Moreover, when $J$ is integrable and $D \subset (M, J)$ is equivalent to a strictly linearly convex domain in $\bC^n$, the family $\cF^{(x_o)}$ of stationary disks 
centered at a fixed  $x_o \in D$,  determines   a smooth  foliation of $D \setminus \{x_o\}$ with several important   properties (\cite{Le, Le1}): 
\begin{itemize}
\item[a)] it  determines a  natural diffeomorphism  $\exp: \overline B^n \longrightarrow \overline D$, which generalizes   the usual Riemann map between $\overline \D$ and any other smoothly bounded domain of $\bC$; 
\item[b)] it consists of     disks that are extremal  for  the Kobayashi metric of $D$; 
\item[c)] it can be used  to determine a  Green pluripotential for $D$ with  pole in $x_o$, 
i.e.  a  plurisubharmonic function that solves  the classical complex Monge-Amp\`ere equation and has   a logarithmic pole at $x_o$. 
\end{itemize}
When $J$ is not integrable, there are several cases,  in which  the family $\cF^{(x_o)}$ of stationary disks,  
centered at a fixed  $x_o \in D$, gives
a smooth  foliation for $D \setminus \{x_o\}$ (e.g. when  $(D, J)$ is  equivalent with a strictly linearly convex domain of $\bC^n$ and $J$  is  a small deformation of $\Jst$). In these cases,  it is still true that  $\cF^{(x_o)}$  determines a $J$-biholomorphically invariant, generalized Riemann map  $\exp: \overline B^n \longrightarrow \overline D$. But  in general,  {\it it is no longer true that the disks of $\cF^{(x_o)}$  are extremal disks for the Kobayashi metric nor that they can be used to  solve  the  almost complex Monge-Amp\`ere equation}. 
Here,  by   ``almost-complex Monge-Amp\`ere equation'' we mean  the differential equation that characterizes the  maximal $J$-plurisubharmonic functions of class $\cC^2$  of an almost complex, strongly pseudoconvex domain. By comparison with usual complex Monge-Amp\`ere equation, it can be considered as a very appropriate analogue  in almost complex settings. 
\par
\smallskip
There are many reasons which justify these phenomena. In case of non-integrable complex structures, the great abundance of $J$-holomorphic curves,  which gives an advantage in many 
geometrical considerations, turns into  a drawback  in considering  objects as the 
Kobayashi metric,  which reveals to be  a weaker and more elusive invariant. 
In particular, it is natural to expect that  the notions of stationary and extremal disks,  which involve    fine (and different!) properties,  become equivalent only when the    almost complex structure satisfies appropriate restrictions. In fact, in \cite{GJ} it is shown that, in general,  these notions are different.  \par
As for the construction of Green pluripotentials,  additional difficulties  emerge. In fact, if one  has the integrable setting in mind, the behavior  of  plurisubharmonic functions    in the non-integrable case is quite unexpected. For instance, 
 there are arbitrarily small deformations $J$ of  the standard complex structure,  w.r.t. which 
the logarithm of squared norm $\log | \cdot |^2$ of $\bC^n$ is  {\it not}  $J$-plurisubharmonic. 
  Furthermore,  the kernel distribution   of an almost-complex Monge-Amp\`ere operator, even if  appropriate non-degenericity conditions is assumed, is usually {\it neither integrable, nor $J$-invariant}.  Since all this is  in clear contrast with the classical  setting,  it  cannot be expected  that for completely arbitrary  non-integrable structures one can  reproduce
the whole  pattern  of  fruitful properties, which relate regular solutions of  complex Monge-Amp\`ere equations and Monge-Amp\`ere foliations.
\par
\smallskip
In this paper,   with the help of  generalized Riemann maps,  we  determine  ``normal forms''  for almost complex domains $(D, J)$ with singular foliations  by stationary disks. We use  such normal form to  construct  examples and determine intrinsic conditions,  under which the disks of the foliations   are   extremal disks  for the Kobayashi metric or give solutions to  the almost complex Monge-Amp\`ere equation.
In fact,  we are able to determine sufficient conditions on the almost complex structure,  which ensure the existence of almost complex Green pluripotential and  the equality between the two notions of stationary disks and of extremal disks. It is interesting to note that the class of such structures (called {\it nice} or {\it very nice} almost complex structures) is very large in many regards,  in fact determined by a finite set of conditions (it is finite-codimensional) in an infinite dimensional space. We hope  that such notions will be fruitful also for 
other questions in almost complex analysis and geometry.
\par
\smallskip
The paper is organized as follows. After a preliminary section, in \S 3 we introduce the notion of 
almost complex domains of circular type {\it in normal form}. They  are pairs $(B^n, J)$, formed by the unit ball $B^n \subset \bC^n$ and an almost complex structure $J$,  which satisfies  conditions 
 that guarantee that  any radial disk through the origin is stationary  for $(B^n, J)$. Since any almost complex domain, admitting a  singular  foliation by stationary disks,   is biholomorphic to a domain   in normal form, 
any problem on such foliations  can be reduced to questions on the radial disks of   normal forms. In \S 4, we study conditions on  $J$,  under which the  radial disks of a normal form $(B^n, J)$ are  extremal. In \S 5, we define the {\it almost complex Monge-Amp\`ere equation}, we prove  that it  characterizes   maximal $\cC^2$ plurisubharmonic functions and we determine conditions  on  normal forms $(B^n, J)$, under which the stationary foliation by radial disks determines a Green pluripotential.   \par
\bigskip
\noindent{\it Notation.} The standard complex structure of $\bC^n$ is denoted by $\Jst$, the unit ball $ \{\ |z| < 1\ \} \subset \bC^n$ is denoted by $B^n$ and, when  $n = 1$, by  $\D = B^1$. \par
For any  $\a> 0$ and  $\epsilon \in ]0,1[$,  a   map $f: \overline \D \longrightarrow M$ into a manifold $M$ is said  {\it of class 
$\cC^{\a, \epsilon}$}
if there are  coordinates $\xi = 
(x^1, \dots, x^N): \cU \longrightarrow \bR^N$ 
on a neighborhood   of $f(\overline \D)$,  such that $\xi \circ f : \overline \D \longrightarrow \bR^n$ is of class $\cC^\a$ on $\D$ and H\"older continuous of class $\cC^\epsilon$ on $\overline
\Delta$.
If $Y= Y_i^j \frac{\partial}{\partial x^j} \otimes d x^i$ is a tensor field of type $(1,1)$ on $\bR^m$ and $\cU$ is a subset  of $\bR^m$, we denote
 $\displaystyle{\| Y \|_{\overline \cU, \cC^k} \=  \sum_{|J| \leq k}\  \sup_{x \in \overline \cU}\ \left|\frac{\partial^{|J|} Y^i_j}{\partial x^J}(x) \right|}$.\par
 \par
\bigskip

\section{Preliminaries}
\label{preliminaries}
\setcounter{equation}{0}
 \subsection{Canonical lifts  of  almost complex structures}
 \label{preliminaries2.2}
 Let  $(M,J)$ be an $n$-dimensional complex manifold with integrable complex structure $J$. In this case, $TM$ and $T^*M$ are naturally endowed with integrable complex structures $\bJ$ and $\wt \bJ$, respectively, corresponding to the  atlases of complex charts 
 $\wh \xi: TM|_{\cU} \longrightarrow T \bC^n \simeq \bC^{2n}$,  $\wt \xi: TM^*|_{\cU} \longrightarrow T^* \bC^n \simeq  \bC^{2n} $,  determined by charts $\xi = (z^i): \cU \subset M \longrightarrow \bC^n$ of the atlas  of the complex manifold structure of  $(M, J)$.  \par
 When $(M, J)$  is  an almost complex manifold, there is no canonical atlas of complex charts on  $M$ and the   above construction is meaningless.  Nevertheless,  {\it there  are   natural almost complex structures $\bJ$ and $\wt \bJ$ on $TM$ and $T^*M$, respectively, also in this more general case} (see \cite{IY}, \S I.5 and \S VII.7). Using coordinates they are defined as follows .  For a 
given a  system of real coordinates
  $\xi = (x^1, \dots, x^{2n}): \cU \subset M \longrightarrow \bR^{2n}$, we denote by  
 \beq \label{coordinates}
\wh \xi = (x^1, \dots, x^{2n}, q^1, \dots, q^{2n}): \pi^{-1}(\cU) \subset TM \longrightarrow \bR^{4n}\ ,\eeq
\beq \label{coordinates1}
\wt \xi = (x^1, \dots, x^{2n}, p_1, \dots, p_{2n}): \wt \pi^{-1}(\cU) \subset T^*M \longrightarrow \bR^{4n}\ ,\eeq
 the associated coordinates on $TM|_{\cU}$ and $T^*M|_{\cU}$,  determined by the components $q^i$ of vectors $v = q^i \frac{\partial}{\partial x^i}$
  and  the components $p_j$ of the covectors $\a = p_j d x^j$.  If   $J^i_j = J^i_j(x)$  denote the  components of
$J = J^i_j \frac{\partial}{\partial x^i} \otimes dx^j$,  the almost complex structures  $\bJ$ and $\wt \bJ$ are defined by  the expressions
\beq  \bJ= J_{i}^a \frac{\partial}{\partial x^a} \otimes dx^i + 
J_{i}^a   \frac{\partial}{\partial q^a}\otimes d q^i  +  
 q^b J^a_{i,b} \frac{\partial}{\partial q^a} \otimes dx^i
\ ,\eeq
$$ \wt \bJ= J_{i}^a \frac{\partial}{\partial x^a} \otimes dx^i + 
J_{i}^a   \frac{\partial}{\partial p_i}\otimes d p_a  +  \phantom{aaaaaaaaaaaaaaaaaaaaaaaaaaaa}
 $$
\beq \label{bJ}  \phantom{aaaaaaaa} + 
\frac 1 2  p_a\left( - J^a_{i,j}  + J^a_{j,i} + J^a_\ell \left(J^\ell_{i,m} J^m_j - J^\ell_{j,m} J^m_i \right)\right)
\frac{\partial}{ \partial p_j}\otimes d x^i\ .\eeq
These tensor fields can be checked to be independent  on the  chart $(x^i)$ and: 
\begin{itemize}
\item[i)] the standard projections $\pi: T^*M \longrightarrow M$,  $\wt \pi: T^*M \longrightarrow M$ are  $(\bJ,J)$-holomorphic and $(\wt \bJ, J)$-holomorphic, respectively;
\item[ii)]   given a $(J,J')$-biholomorphism $f: (M, J)\longrightarrow (N,J')$ between almost 
complex manifolds, the tangent and cotangent  maps 
$$ f_*: T M \longrightarrow T N\qquad \text{and} \qquad   f^*: T^*N \longrightarrow T^*M$$
are $(\bJ, \bJ')$- and  $(\wt \bJ',\wt \bJ)$-holomorphic, respectively;
\item[iii)] when $J$ is integrable,   $\bJ$ and $\wt \bJ$  coincide with above described integrable complex structures of $TM$ and $T^*M$, respectively (in the integrable case, all derivatives $J^a_{i,j}$ are  $0$ in  holomorphic coordinates). 
\end{itemize}
\par
\smallskip
We call  $\bJ$, $\wt \bJ$ 
 {\it  canonical 
lifts of $J$ on $TM$ and $T^*M$}. \par
 \par
 \medskip 
 \subsection{Blow-ups  of almost complex manifolds} 
 Given a point $x_o$ of an almost complex manifold $(M, J)$, we call   {\it blow-up at $x_o$} the topological manifold $\wt M$    obtained as follows (\cite{PS1}).  \par
  Consider a system of complex coordinates $\xi = (z^1, \dots, z^n): \cV \longrightarrow \cU \subset \bC^n$ on a  neighborhood $\cV$ of $x_o$,  with  $\xi(x_o) = 0$ and which maps $J|_{x_o}$ into the standard complex structure $\Jst|_0$ of 
$T_0 \bC^n \simeq \bC^{2n}$. The manifold  $\wt M$ is obtained   by   gluing  $M \setminus \{x_o\}$ with the blow up 
 $\wt \cU$ at $0$ of  $\cU \subset \bC^n = \bR^{2n} $,  identifying  $\cV \setminus \{0\}$ with $\cU \setminus \{0\}$  by means of  the map $\xi = (z^i)$.  As it was remarked in \cite{PS1},  the smooth manifold structure of  $\wt M$ does not depend on the choice of the coordinates $\xi = (z^i)$.  Hence, this manifold $\wt M$ can be considered as canonically associated with $(M, J)$ and $x_o$.  \par
Since $M \setminus \{0\} \equiv \wt M \setminus \pi^{-1}(x_o)$ and  $ \pi^{-1}(\cV) \equiv \pi^{-1}(\wt \cU) \subset \wt B^n$, we may consider  the tensor field $J$ of type $(1,1)$ on $\wt M$, with   $J_x: T_x\wt M \to T_x \wt M$  equal to  the almost complex structure of $M$ for  any point $x \in \wt M \setminus \pi^{-1}(x_o)$ 
and to  the standard complex structure   of $\wt B^n$ for   any $x \in \pi^{-1}(x_o) \equiv \pi^{-1}(0) \subset \wt B^n$.  Such  tensor field is obviously smooth on $M \setminus \pi^{-1}(x_o)$ and, in our discussions, there will be  no need to know whether  it is  smooth   also on $\pi^{-1}(x_o)$.  It is however possible to check that it  is  in fact smooth at all points.    \par

 \par
\bigskip
\section{Normal forms of almost complex domains of circular type}
  \setcounter{equation}{0}
    In what follows, $D \subset M$ denotes a domain in a $2n$-dimensional almost complex manifold $(M, J)$  with smooth boundary $\Gamma = \partial D$.\par
\medskip
\subsection{Almost complex domains of circular type }
Let    $\cN$ be the conormal bundle of $\Gamma = \partial D$, i.e. the subset of $T^*M|_{\G}$
$$\cN = \{\ \b \in T^*_xM\ ,\ x \in \G\  :\ \ker \b \subset T_x \G\ \}\ .$$
 We recall that, given $\a\geq 1$ and $\varepsilon> 0$, a    {\it $\cC^{\a,\varepsilon}$-stationary disk  of $D$\/}  is a map  $f: \overline \Delta \longrightarrow M$ such that
\begin{itemize}
\item[i)] $f|_{\Delta}$ is  a $J$-holomorphic embedding and  $f(\partial \Delta) \subset \partial D$; 
\item[ii)]  there exists  a  $\wt \bJ$-holomorphic  map $\wt f :  \overline \Delta \longrightarrow T^*M$ with $\pi \circ \wt f = f$,  so that 
\beq \label{stationarycondition} \zeta^{-1}\cdot \wt f(\zeta) \in \cN\setminus\{\text{zero section}\}\ \  \text{for any} \ 
\zeta \in \partial \Delta\eeq
 and  $ \wt \xi \circ \wt f \in \cC^{\a, \varepsilon}(\overline \D, \bC^{2n})$ for some complex coordinates $\wt \xi = (z^i, w_j)$  around  $\wt f(\overline \D)$. 
 \end{itemize}
 In \eqref{stationarycondition}  ``\ $ \cdot$\ "  denotes   
the usual $\bC$-action  on $T^*M$, i.e. the action 
(\footnote{We follow Besse's convention on signs, for which $J^* \a(v) = - \a(Jv)$ (\cite{Bs}; see  also \S \ref{plurisubharmonic}). Due to this,  on $\bC^n$ we have $\Jst^* dx^j = dy^j$,  $d z^j = dx^j + i \Jst^* dx^j$ and $ i d z^j = -\Jst^* dz^j$.})
\beq\label{product} \z\cdot \alpha \= \Re(\z)\alpha  -  \Im(\z) J^*\alpha\quad \text{for any}\ \alpha \in T^* M, \ \ \z \in \bC\ \ .\eeq
 If $f$ is stationary, the maps $\wt f$ satisfying  (ii) are called  {\it stationary lifts of $f$}.\par
 \medskip
 In this paper we are concerned with domains $D$ in an almost complex manifold $(M, J)$, admitting singular foliations by stationary disks with the same properties of the singular foliations  by Kobayashi extremal disks of  the domains of circular type in $\bC^n$.  Here is the definition of such domains.\par
\begin{definition} \cite{PS1} \label{foliation} For any point     $x_o$ of an almost complex domain $D \subset (M,J)$,  we 
 denote by $\cF^{(x_o)}$ the family  of stationary disks of $D$ with  $f(0) = x_o$. We say the   $\cF^{(x_o)}$  is a {\it  foliation of circular type\/}  if:
\begin{itemize}
\item[i)] for any $v \in T_{x_o} D$, 
 there exists a unique disk $f^{(v)} \in \cF^{(x_o)}$
with
$f^{(v)}_*\left( \left.\frac{\partial}{\partial x}\right|_0 \right) = \mu \cdot v$ for some $0 \neq  \mu \in \bR$;
\item[ii)]  for a fixed identification  $(T_{x_o} D, J_{x_o}) \simeq (\bC^n, \Jst)$,  
 the map from the blow up  $\wt B^n$ at  $0$ of   $B^n$ to  the blow up $\wt D$ at $x_o$ of  $D$
\beq \label{exp} \exp: \wt B^n  \subset \wt \bC^n \longrightarrow  \wt D \ \ , \qquad  \ \exp(v, [v]) \= \wt{f^{(v)}}(|v|)\eeq
is smooth, extends smoothly up to the boundary and determines a diffeomorphism  between the boundaries  $\exp|_{\partial B^n}: \partial B^n \longrightarrow \partial D$.
 \end{itemize}
The point   $x_o$  is called  {\it center of the foliation\/} and     $\exp: \wt B^n \longrightarrow \wt D$  is called {\it (generalized) Riemann map of $(D, x_o)$}.  Any domain $D \subset (M, J)$ admitting a foliation of circular type is called  {\it  almost complex domain  of  circular type}.\end{definition}
There  exists of a wide class of almost complex domains of circular type.  In fact,   any bounded, strictly linearly convex domains  $\wt D \subset \bC^n$,   with smooth boundary and endowed with a small deformation  $J$ of the standard complex structure $\Jst$, admits singular  foliations of circular type made of  stationary disks   (\cite{PS1}, Thm. 4.1). For  such domains,  Gaussier and Joo proved  in \cite{GJ} the same existence result  for singular foliations,  made of the  so-called {\it $J$-stationary disks}  (see   later  for the definition).\par 
\medskip
On the other hand,  the following fact  is well-known (see  e.g. \cite{DS}). 
 \begin{lem} \label{small} Let $(M,J)$ be an almost complex manifold  and $x_o\in M$. For any integer   $k \geq 0$ and  $\varepsilon > 0$, 
  there exists  a neighborhood $\cU$ of $x_o$,   such that  $(\cU, J)$ is   $(J, J')$-biholomorphic to  $(B^n , J')$ for some  almost complex structure   $J'$ on a neighborhood of   $\overline B^n \subset \bC^n$ with  $\| J' - \Jst\|_{\overline{B^n},\ \cC^k} < \varepsilon$. 
 \end{lem}
This and previous remarks imply that  any point $x_o$ of an almost complex domain $(M, J)$ admits  a neighborhood $\cU$  containing a  domain $D \subset \subset \cU$ of circular type with center $x_o$.  \par
\medskip
\subsection{Normal forms} 
\label{domainsofcirculartype}
By definitions, 
for any almost complex domain $(D, J)$  of circular type with center $x_o$, the map  $\exp: \wt B^n \subset \wt \bC^n \longrightarrow \wt D$ is a biholomorphism 
between $(\wt D, J)$ and    $(\wt B^n, \wt J)$,  where   $\wt J\= \exp^{-1}_*(J)$. \par
If we denote by $\pi:\wt B^n \longrightarrow B^n$,   $\pi':\wt D \longrightarrow D$ the natural blow down maps and by $J' = \pi_*(\wt J)$ the projected almost complex structure on $B^n\setminus \{0\}$,  we have that the map $E = \pi \circ \exp \circ \pi'{}^{-1}: D\setminus \{0\}) \longrightarrow B^n\setminus \{0\}$ is a $(J,J')$-biholomorphism.  In general,  the tensor field $J'$ does not extend   smoothly  at $0 \in B^n$. Nonetheless  such singularity  is  ``removable''  in  the following sense. \par
First of all,  we remark that  the map $E$ extends uniquely to a {\it homeomorphism}  $E: D\longrightarrow B^n$ by setting  $E(x_o) = 0$. So,   we  may  consider the atlas  $\cA$ on $B^n$,  formed by the charts  $\h = \xi \circ E^{-1}$  determined by  charts $\xi: \cU \subset D \to \bR^{2n}$  of the manifold structure of $D$.  {\it Such atlas defines a smooth manifold structure on $B^n$,  which coincides with the standard  one on $B^n \setminus \{0\}$, but contains   charts around  $0$  that in general are non-standard}.  By construction,  the components of $J'$ in the charts of $\cA$ extend smoothly    at  $0 $. \par
\medskip
We call the  pair   $(B^n, J')$   {\it normal form
of $(D, J)$   determined by $\cF^{(x_o)}$}. If   we endow $B^n$ with the atlas $\cA$,  by construction  $(B^n, J')$ is an almost complex domain of circular type with center $0$ and foliation $\cF^{(0)}$ given by the {\it straight disks} 
\beq\label{straightdisk}Ê f^{(v)} : \D\longrightarrow B^n\ ,\qquad f^{(v)}(\z) = \z \cdot v\ ,\qquad v \in \bC^n\ .\eeq
\par
\medskip
Now, let us  give an  intrinsic characterization of the  almost complex structures $J$ on $B^n\setminus \{0\}$ that correspond to  normal forms of domains of circular type. For this purpose, we need  some new notation and the notion  of ``almost $L$-complex structures''. 
 \par
 \smallskip
 Let  $Z  \= \Re\left(z^i \frac{\partial}{\partial z^i}\right)$ and denote by   $\cZ$ the $\Jst$-invariant distribution on $B^n\setminus\{0\} \subset \bC^n$ defined  by $\cZ_z \= < Z_z, \Jst Z_z> $ at any $z \neq 0$. We recall  that  
 \beq\cZ_z   = \ker \left.d d^c_{\operatorname{st}} \log \tau_o\right|_x\ ,\quad\text{where} \ \  \t_o(z) \= |z|^2\ , \ \ \dcst \= \Jst^* \circ d \circ \Jst^*\label{2.14}\eeq 
 and that
 $d\dcst\tau_o(Z, X)  = X(\tau_o)$ for any $X \in T B^n \setminus \{0\}$. One can  check that $\cZ$ is integrable and that its  integral leaves are   the (images of the) disks \eqref{straightdisk}.  \par
 Consider the blow up $\wt B^n$ of $B^n$ at $0$, the standard identification of $\wt B^n$ with an open subset of the tautological bundle $\pi: E \longrightarrow \bC P^{n-1}$ and the coordinates $\xi : \cU\subset E\longrightarrow \bC^n$,  $\cU \=  \{\ ([v], v) : \ v^n \neq   0\ \} $, defined by
\beq \label{interestingcoord} 
\xi^{-1}(z^0, \dots, z^{n-1}) \= \phantom{aaaaaaaaaaaaaaaaaaaaaaa}$$
$$\= \left(\!\left[z^1 :\dots : z^{n-1} : \sqrt{1- \sum_{i = 1}^{n-1} |z^i|^2}\right]; z^0\!\!\cdot\!\! \left(z^1, \dots, z^{n-1}, \sqrt{1- \sum_{i = 1}^{n-1} |z^i|^2}\right)\!\!\right)\!\!.\eeq
In these coordinates,
$\cZ^\bC_z = \Span_{\bC}\left\{\ \left.\frac{\partial}{\partial z^0}\right|_z, \left.\frac{\partial}{\partial \overline{z^0}}\right|_z\  \right \}$ 
and $J \left(\left.\frac{\partial}{\partial z^0}\right|_z\right) = i \left.\frac{\partial}{\partial z^0}\right|_z$.
Now, let us use capital letters  $A, B, C, \dots$ for indices  that might be indifferently of the form $a$ or $\overline a$,  so that  we may denote the complex coordinates and their conjugates by  $(z^A)= (z^a, z^{\overline a} \= \overline {z^a})$. Let also denote by  $(p_A) = (p_a, p_{\overline a} \= \overline{p_a})$ the  complex components of real  1-forms $\o = p_a dz^a + \overline p_{a} d z^{\overline a} \in T^* B^n$.  Using these conventions,  the canonical lift $\wt \bJ$ of an almost complex structure $J$  on $B^n \setminus \{0\}$ is of  the form 
$$\wt \bJ = J^B_A \left(\frac{\partial}{\partial z^B} \otimes d z^A + \frac{\partial}{\partial p_A} \otimes dp_B\right) + \phantom{aaaaaaaaaaaaaaaaaaaaaaaaaaa}$$
$$ \phantom{aaaaa}+ \frac{1}{2} p_C \left( - J^C_{A, B} + J^C_{B, A} + J^C_L\left( J^L_{A, M} J^M_B - J^L_{B, M} J^M_A\right)\right) \frac{\partial}{\partial p_B} \otimes dz^A\ ,$$ 
where $J^A_B$ are the components of $J$ w.r.t. the complex vector fields $\left(\frac{\partial}{\partial z^A}\right)$.  One way to recover such  formula is, for instance,  to look at the expression of $\wt \bJ$ in terms of Nijenhuis  tensor and tautological form of  $T^*M$ (\cite{IY, Sp}) and write  all terms using the  complex  coordinates $(z^A, p_B)$. Again, we remark that,  when  $J$ is integrable and $(z^1, \dots, z^n)$ are holomorphic coordinates,  $\wt \bJ = J^B_A \left(\frac{\partial}{\partial z^B} \otimes d z^A + \frac{\partial}{\partial p_A} \otimes dp_B\right) $\par
\medskip
We can now introduce the ``almost $L$-complex structures'': as we will shortly see, a pair $(B^n, J)$ is a domain of circular type in normal form if and only it   $J$ is an almost complex structure of this kind (Theorem \ref{chara}). We will also see  that these structures are characterized by a finite collection of equations in the space of almost complex structures (Proposition  \ref{usefulrem}). \par
 
\begin{definition} \label{Lstructure} We call {\it almost $L$-complex structure}  any  almost complex structure  $J$   on $B^n\setminus \{0\}$, smoothly extendible at $\partial B^n$, such that: 
\begin{itemize}
\item[i)]   $\cZ$ is $J$-stable  and   $J|_{\cZ} = \Jst|_{\cZ}$;  
\item[ii)] for any $v \in S^{2n-1}$,  the following differential problem on $2n-2$ $\bC$-valued maps $g_\a, g_{\bar \a}: \overline \D \longrightarrow \bC$ of class $\cC^{\a, \epsilon}$  is solvable (in \eqref{diffsyst},  $A$, $B$ denote indices of the form $\a$, $\bar \a$, $\b$, $\bar \b$, respectively, with $1 \leq \a, \b \leq n-1$)
\beq \label{diffsyst} \left\{ \begin{array}{lr}   \left(\d_A^B 
- i (\left.J_A^B \right|_{\z \cdot v}) \right) g_{B, \overline\z} + & \ 
\\
\phantom{aaaaaaaaaaa}
+   \left( -  \frac{i}{2} \left.\left( J^B_{A, \overline 0} + i J^B_L J^L_{A, \overline 0}\right) \right|_{\z \cdot v}\right)  g_B + & \ \\ \phantom{aaaa}
- \left( \frac{i}{2}  \left.\left( J^0_{A, \overline 0} + i J^0_L J^L_{A, \overline 0} + J^{\bar 0}_{A, \overline 0} + i J^{\bar 0}_L J^L_{A, \overline 0} \right) \right|_{\z \cdot v}\right)  = 0& \text{when}\ \z \in \D\ ,\\
 \ & \ \\
 \left((\Re \z) \d_A^B -  (\Im \z)   \left. J^B_A\right|_{\z \cdot v} \right)g_B   -  & \ \\
\phantom{aaaaaaaaaaaaaaaaataaa} - (\Im \z) \left.\left(  J^0_A + J^{\bar 0}_A\right) \right|_{\z \cdot v}= 0  & \!\!\! \text{when}\ \z \in \partial \D \ ,\end{array}\right.\eeq
where $(J^B_A)$ are the components of $J$ in coordinates  of the form \eqref{interestingcoord};
\item[iii)] there exists a homeomorphism $\xi: \cU \longrightarrow \cV$ between neighborhoods  of  $0 \in \bC^n$, which is   $\cC^{\infty}$ on $\cU \setminus \{0\}$ and such that  $\xi_* (J)|_{\cV \setminus \{\xi(0)\}}$ extends smoothly  at $0$;  in particular,   $J$ admits a smooth extension at $0$ if $B^n$ is endowed with a (non-standard) atlas containing $\xi$; 
\item[iv)] the blow-up $\wt B^n$ of $B^n$, determined by $J$ and  the non-standard smooth manifold  structure described in (iii),   is diffeomorphic to the usual blow-up of $B^n$ 
determined by $\Jst$. 
\end{itemize}
 Unless explicitly stated,  for any almost $L$-complex structure, we will always assume $B^n$ endowed with the  smooth manifold  structure described   in (iii). 
\end{definition}
It is useful to remark that condition (ii) of the above definition is satisfied by a very large class of almost complex structures. Moreover, in the integrable case,  (ii) is automatically satisfied by the complex structures of the normal forms of domains of circular types (see \cite{PS}). \par
\begin{prop} \label{usefulrem} Let $J$ be  an almost complex structure  $J$   on $B^n\setminus \{0\}$ such that,  for any $v \in S^{2n-1}$ and $\z \in \partial \D$, the matrices
 $$\gC = \left[\left(\d_A^B 
- i (\left.J_A^B \right|_{\z \cdot v}) \right) \right]\ ,\quad \gD = \left[ \left((\Re \z) \d_A^B -  (\Im \z)   \left. J^B_A\right|_{\z \cdot v} \right)\right]$$
are  invertible.  Let also $\bF \in  \cC^{\a-1, \epsilon}(\overline \D, \bC^{2n-2})$ and $\bG \in \cC^\epsilon (\partial \D, \bR^{4n-4})$ defined by 
$$\bF(\z) = \gC^{-1}(\z) \cdot \left( \frac{i}{2}  \left.\left( J^0_{A, \overline 0} + i J^0_L J^L_{A, \overline 0} + J^{\bar 0}_{A, \overline 0} + i J^{\bar 0}_L J^L_{A, \overline 0} \right) \right|_{\z \cdot v}\right)\ ,$$
$$\bG(\z) = (\Im \z) \left.\left(  J^0_A + J^{\bar 0}_A\right) \right|_{\z \cdot v}\ .$$
Then,  \eqref{diffsyst} is solvable if and only if
 $(\bF, \bG)$ is in the (finite codimensional) range of the Fredholm operator  described in formula \eqref{fredholm} below.
 \end{prop} 
 \begin{pf} The system \eqref{diffsyst} is equivalent to the  ``generalized Riemann-Hilbert problem'' on maps $g: \overline \D \longrightarrow \bC^{2n-2}$
\beq \label{diffsystbis} \left\{ \begin{array}{lr}  g,_{ \bar\z} + (\gC^{-1}  \bA)\cdot g = \bF
& \text{on}\  \D\ ,\\
 \ & \ \\
\gD\cdot  g = \bG  & \!\!\! \text{on}\ \  \partial \D \ ,\end{array}\right.\eeq
where  $\bA: \D \to M_{2n - 2}(\bC)$ is 
$\bA(\z) \=  \left[ -  \frac{i}{2} \left.\left( J^B_{A, \overline 0} + i J^B_L J^L_{A, \overline 0}\right)\right|_{\z \cdot v}\right]$.
 The operator 
 \beq \label{fredholm}ÊR: \cC^{\a, \epsilon}(\overline \D, \bC^{2n-2}) \longrightarrow \cC^{\a-1, \epsilon}(\overline \D, \bC^{2n-2}) \times \cC^\epsilon (\partial \D, \bR^{4n-4})\phantom{aaaaaaaaaaaa} $$
 $$ \phantom{aaaaaaaaaaaaa} R(h) \=\left( \frac{\partial h }{\partial \z} + \bA \cdot h ,  \left(\Re  \left(\gD\cdot  h\right), \Im \left(\gD\cdot  h\right) \right)\right)\eeq
 is known to be  Fredholm if and only if $\gD$ is invertible at all points (\cite{Wen} \S 3.2 and \cite{MP} \S VII.3). The claim follows immediately.
 \end{pf}
 \begin{rem} \label{usefulrembis} Notice that if the components $J_A^0$, $J^{\bar 0}_A$, appearing in \eqref{diffsyst}, are identically equal to $0$ along the considered disk,  the system    \eqref{diffsyst} always admits the trivial solutions $g_\a \equiv 0 \equiv g_{\bar \a}$, regardless on the invertibility of $\gC$ and $\gD$. This fact turns out to be quite useful to produce examples. 
 \end{rem}
The interest for almost $L$-complex structure is motivated by the following. \par
\begin{lem} \label{charanorm}  If  $J$ is an almost $L$-complex structure on $B^n$, any
  straight disk  through $0$ is stationary w.r.t. $J$.
\end{lem}
\begin{pf} Using coordinates \eqref{interestingcoord}, since $J$ satisfies (i) of Definition \ref{Lstructure} , then 
\beq \label{firstcond} J^A_0 = i \d^A_0\ ,\quad J^A_{\overline 0} = - i \d^A_{\overline 0}\ ,\quad J^A_{0,B} = 0\ ,\quad J^A_{\overline 0, B} = 0\ .\eeq
Given $v = (v^1, \dots, v^{n-1}) \in \bC^{n-1}$, consider the straight disk $f:\D \longrightarrow B^n$, with tangent direction at $0$ given by 
  $\left[v^1: \dots: v^{n-1}:\sqrt{1- \sum_{i = 1}^{n-1} |v^i|^2}\right]$, i.e.  
$$f (\z) \= \xi^{-1}(\z, v^1, \dots, v^{n-1})\qquad \text{for any}\ \  \z \neq 0\ .$$
Recall that  $f$ is stationary if and only if  there is
$$\wt f: \overline \D \longrightarrow T^* \wt B^n\ ,\qquad \wt f(\z) = (\z, v^1,  \dots, v^{n-1}; g_A(\z)_{A = 0, \bar 0, 1, \bar 1, \dots, n-1, \overline{n-1}}) $$
such that: a)  $\z^{-1} \cdot \wt f(\z)$ is in the conormal bundle $\cN \setminus \{\text{zero section}\}$ of $\partial B^n$  for any $\z \in \partial \D$; b)   $\wt f$ is $\wt \bJ$-holomorphic, i.e. 
$ 
\wt f_*\left(\Jst \frac{\partial}{\partial \overline \z}\right) = \wt \bJ\left(\wt f_*\left( \frac{\partial}{\partial \overline \z}\right)\right)$.\par
We claim that a  map $\wt f$  of the above form and satisfying (a) and (b),    exists.  In fact, by   \eqref{firstcond},   condition (b)  is equivalent to (here indices $A$, $B$ might assume any value, $0$, $\bar 0$ included)
  \beq \label{secondcond} - i g_{A, \overline \z} -  (\left.J_A^B \right|_f)  g_{B, \overline\z} - 
  \frac{1}{2} g_B \left(\left.\left( J^B_{A, \overline 0} + i J^B_L J^L_{A, \overline 0}\right) \right|_f\right) = 0\ .\eeq
By  \eqref{firstcond}, in case  $A = 0$ equation   \eqref{secondcond}  reduces to 
$ 2 i g_{0, \overline \z}  = 0$, i.e.  to the requirement  of holomorphicity for $g_0: \D \longrightarrow \bC$.
On the other hand,  since the conormal bundle  $\cN$   is generated at any point by  $\o = z^0 d \overline z^0 + \overline z^0 d z^0$,  (a)  is equivalent to the existence of a continuous  $\l: \partial \D \longrightarrow \bR \setminus \{0\}$ such that  
\beq \label{thirdcond} \z^{-1} g_0(\z)  =  \overline{ \z} \l(\z)\ \ \text{for any} \ \ \z \in \partial \D \ ,  \eeq
(and hence  $g_0 \equiv 1 \equiv g_{\overline 0}$)  and  to the requirement that  $g_\a|_{\partial \D}$,  $g_{\bar \a}|_{\partial \D}$, $\a \neq 0$,  satisfy the boundary conditions of \eqref{diffsyst}. Therefore,  inserting  $g_0 = g_{\bar 0} = 1$ into \eqref{secondcond},   by condition (ii) of Definition \ref{Lstructure},  we conclude that there  always exists a stationary lift    $\wt f(\z) = (\z, v^1,  \dots, v^{n-1}; g_A(\z)) $ for  $f$. 
\end{pf}
By Lemma \ref{charanorm}, if $J$ is an almost $L$-complex structure, $(B^n, J)$ is an almost complex domain of circular type that coincides with its normal form. Conversely, one can directly  check that if $(B^n, J)$ is a domain of circular type in normal form,  then $J$ satisfies all  conditions of Definition \ref{Lstructure}. We have therefore the following intrinsic characterizations of normal forms. 
\begin{theo} \label{chara} A pair $(B^n, J)$  is   a domain of circular type in normal form if and only if $J$ is an almost $L$-complex structure. 
\end{theo}
\bigskip
\subsection{Deformation tensors  of normal forms}
Consider now the  $\Jst$-invariant distribution 
  $\cH$   on $B^n \setminus \{0\}$
 defined by 
\beq \cH_z  \= \{\ X \in T_xM\ :\ d\dcst\tau_o(Z, X) = d\dcst\tau_o(\Jst Z, X) = 0 \ \}\ .\label{2.12}\eeq
 One can directly check that   $T_z M = \cZ_z \oplus \cH_z$ and that    $\cH_z$  coincides with  the holomorphic tangent  space to the sphere $S_{|z|} = \{\ w\ : |w| =|z|\ \}$. In particular, for any $0 < c < 1$, 
 {\it the pair  $(\cH|_{S_c}, \Jst)$ is   the CR structure of the  sphere  $S_c  = \{\ \tau_o = c\ \}$\/}.
It is known that  the distributions $\cZ$ and $\cH$  extend smoothly  on 
the blow up $\widetilde{B}$ and that also such extensions are  $\Jst$-invariant  (see e.g. \cite{Pt, Pt1, PS}). \par
\smallskip
Recall that any complex structure $J_z$ on a tangent space $T_z B^n$ is uniquely determined by  its $-i$-eigenspaces $(T_z B^n)^{01}_J$ in $T^\bC_z B^n$. If  $\cZ^\bC_z = \cZ^{10}_z + \cZ^{01}_z$ and $\cH^{\bC}_z = \cH^{10}_z + \cH^{01}_z$ are decompositions into $\Jst$-eigenspaces,  a generic complex structure $J_z:T_z B^n \longrightarrow T_z B^n$ is completely determined by the   tensors 
$$\phi^\cZ_z \in\Hom( \cZ^{01}_z,  \cZ^{10}_z)\ ,\quad  \phi^\cH_z \in\Hom( \cH^{01}_z,  \cH^{10}_z)\ , \quad 
\phi^{\cZ,\cH}_z \in\Hom( \cZ^{01}_z,  \cH^{10}_z)\ ,$$
$$ \phi^{\cH,\cZ}_z \in\Hom( \cH^{01}_z,  \cZ^{10}_z)\ ,$$
 which determine the $-i$-eigenspace $ (T_z B^n)^{01}_J $ as the complex subspace 
\beq \label{decomp} (T_z B^n)^{01}_J  =\!\! \left(\phantom{A^{A^A}}\!\!\!\!\!\!\!\!\!\!\!\!\cZ^{01} + \phi^\cZ_z(\cZ^{01}) + \phi^{\cZ,\cH}_z(\cZ^{01})\right) +
 \left(\phantom{A^{A^A}}\!\!\!\!\!\!\!\!\!\!\!\!\cH^{01} + \phi^\cH_z(\cH^{01}) + \phi^{\cH,\cZ}_z(\cH^{01})\right)\!\!. \eeq
We call {\it deformation tensor of $J$ w.r.t. to $\Jst$} the tensor field $\phi \in (T^{01 *} \otimes T^{10})(B^n \setminus \{0\})$,  defined at any point  by  $\phi_z \= \phi^\cZ_z  + \phi^{\cZ,\cH}_z  +
 \phi^\cH_z + \phi^{\cH,\cZ}_z$.
From Definition \ref{Lstructure} (i),  it follows immediately  that
\begin{prop} A generic  almost $L$-complex structure $J$ is uniquely determined by  a deformation tensor  of the form  
 \beq \label{deform1} \phi_z =
 \phi^\cH_z + \phi^{\cH,\cZ}_z\qquad \text{\it for any} \ z \in B^n \setminus\{0\}\ .\eeq
and, conversely,  any deformation tensor  as in \eqref{deform1} gives an almost $L$-complex structure,   provided that  the corresponding  $J$ extends to $\partial B^n$ and  $0$ as required in  (ii)  -  (iv)  of
 of Definitions \ref{Lstructure}.
 \end{prop}
 
 \begin{rem}\label{remark} Conditions (iii) - (iv) of Definition \ref{Lstructure} are requirements  that might be hard to check. However,  to construct  examples, it is often sufficient  to observe that,    given a deformation tensor  $\phi_o$  that satisfy those conditions  (e.g.  $\phi_o \equiv 0$), also the  deformation tensors of the form $\phi = \phi_o + \d \phi$,  in which  $\d \phi$  vanishes  identically   on some neighborhood $\cU$ of $0$, satisfy them.\par
 Notice also that if  $J$ is an almost complex structure, determined by a deformation tensor of the form  $\phi = \phi^\cH$, for any given  disk $f(\z)= \z \cdot v$, we may choose coordinates of the form \eqref{interestingcoord}, in which $v = (0, \dots, 0, 1)$ and hence  the components $J_A^0$, $J^{\bar 0}_A$, appearing in \eqref{diffsyst}, are identically equal to $0$ along the disk. In this case,   \eqref{diffsyst} always admits the trivial solutions $g_\a \equiv 0 \equiv g_{\bar \a}$. 
 On the other hand, condition (ii) corresponds to the solvability of  the system that determines stationary lifts and hence it is independent on the choice of the coordinate system. All this implies  that when $\phi = \phi^\cH$, condition (ii)  of Definition \ref{Lstructure} is always automatically satisfied. \end{rem}
 \par \bigskip
\section{Extremal disks and critical foliations}
  \setcounter{equation}{0}
\subsection{Critical and extremal disks } In this section, we recall    the notion  of  ``critical disks'',   recently  introduced by Gaussier and Joo in \cite{GJ}  in their studies on the  extremality w.r.t. the Kobayashi metric of $J$-holomorphic  disks. We have to point out that ``critical disks''   is not the name  used in \cite{GJ}. In that paper, such disks are  called   ``disks vanishing the first order variations''.\par
\smallskip
For this,  we first need to remind of  a few concepts related with the geometry of the   tangent bundle of a manifold $M$. We recall that   the {\it vertical distribution  in $T(TM)$}  is the subbundle   of $T(TM)$ defined by 
 $$T^V (TM) = \bigcup_{(x,v) \in TM} T^V_{(x,v)} M\ \ , \qquad T^V_{(x,v)} M = \ker \pi_*|_{(x,v)}\ .$$
 For any $x \in M$,  let us denote by $(\cdot)^V: T_x M \longrightarrow T^V_{(x,v)} M$ the map 
$\left(\left.w^i \frac{\partial}{\partial x^i}\right|_x\right)^V \= w^i \left.\frac{\partial}{\partial q^i}\right|_{(x,v)}$.  It is possible to check that this map does not depend  on   the choice of   coordinates and that it  determines a natural map  from $TM$ to $T (TM)$ (see \cite{IY}). For any $w \in TM$, the corresponding vector $w^V \in T(TM)$ is called {\it vertical lift of $w$}. \par
\begin{definition}\cite{GS}  \label{variationalfields}  Let $f :  \overline \Delta \longrightarrow M$ be a $\cC^{\a, \epsilon}$,  $J$-holomorphic embedding with   $f(\partial \Delta) \subset \partial D$.  We  call {\it infinitesimal variation of $f$} 
any
 $\bJ$-holomorphic map $W: \overline \D \longrightarrow T M$ of class  $\cC^{\a-1, \epsilon}$ with $\pi \circ W = f$ (here, $\pi: TM \longrightarrow M$ is the natural projection).  An infinitesimal variation $W$  is called {\it attached to $\partial D$ and  with  fixed center}  if 
\begin{itemize}
\item[a)] $\a(W_\z) = 0$ for any $\a \in \cN_{f(\z)}$, $\z \in \partial \D$,
\item[b)] $W|_0 = 0$.
\end{itemize}
It is called {\it with fixed central direction} if in addition it satisfies 
\begin{itemize}
\item[c)] $ \left.W_*\left(\frac{\partial }{\partial \Re \z}\right|_0\right)  \in T^V_{W_0}(TM)$ and it is equal to $ \l \left( f_*\left( \left.\frac{\partial }{\partial \Re \z}\right|_0\right) \right)^V$ for some $\l \in \bR$.
\end{itemize}
The disk   $f$ is called  {\it critical}  if  for any infinitesimal variation $W$,  attached  to  $\partial D$ and with fixed central direction, one has  $ \left.W_*\left(\frac{\partial }{\partial \Re \z}\right|_0\right) = 0$.
\end{definition}
\begin{rem} \label{remark3.2}
The previous definition  is motivated by the following facts. When $f^{(t)}:  \overline\D \longrightarrow M$, $t \in ]-a, a[$,  is a smooth 1-parameter family of $J$-holomorphic disks of class $\cC^{\a, \epsilon}$ with $f^{(0)} = f$, it is simple to check that  
$W \= \left.\frac{d f^{(t)}}{dt}\right|_{t = 0}$ 
is a variational field on $f$. Moreover, if $f^{(t)}$ is such that, for all $t \in ]a,a[$
\beq  f^{(t)}(\partial \D) \subset \partial D\ ,\ \  f^{(t)}(0) = f(0)\ , \ \  f^{(t)}_*\left( \left.\frac{\partial }{\partial \Re \z}\right|_0\right) \in \bR f_*\left( \left.\frac{\partial }{\partial \Re \z}\right|_0\right)\ ,\eeq
then  $W$ satisfies (a) - (c).   On the other hand,  a disk
$f$ is a {\it  locally extremal disk} if  for any  $J$-holomorphic disk $g: \overline \D \longrightarrow M$ of class  $\cC^{\a, \epsilon}$, with  image  contained  in  some neighborhood of $f(\overline \D)$
and such that, for some  $\l \in \bR$, 
$$ g(\partial \Delta) \subset \partial D\ , \quad g(0) = f(0) = x_o\ , \quad g_*\left(\left.\frac{\partial }{\partial \Re \z}\right|_0\right) = \l f_*\left(\left.\frac{\partial }{\partial \Re \z}\right|_0\right)\ ,$$
then  $ \l \leq 1$.  One can directly check that {\it any   locally extremal disk $f$, with $f(\partial \D) \subset \partial D$,  is critical} (see  \cite{GS}, proof of Thm. 4.3).  Conversely, by \S 5 and  Thm. 6.4 of \cite{GJ},   {\it in case $D \subset \bC^n$ is strictly convex  on a neighborhood of $f(\overline \D)$ (in suitable cartesian coordinates) and $J$ is sufficiently close to $\Jst$,  any critical disk  $f$  is locally extremal}.\par
\smallskip
Notice that, when $J$ is close to $\Jst$, the  critical disks  are characterized by  properties that are closely resembles  to those that define the stationary disks.  Disks with such properties are called {\it $J$-stationary disks}  (\cite{GJ}).\end{rem}
\medskip
It is  well-known that,  when  $J$ is integrable, the disks that are stationary coincide with the disks that are critical (see e.g.  \cite{Le,Pa}). This equality  is no longer valid for  generic non-integrable complex structures.   In  \cite{GJ}, counterexamples are given. \par
\smallskip
Next  theorem gives conditions that imply the equality between stationary and critical disks and will be used in the sequel.  The claim and the proof are refinements of a result and arguments given in  \cite{GJ0}.
 In the statement,     $f :  \overline \Delta \longrightarrow M$ is a  $J$-holomorphic embedding,  of class $\cC^{\a, \epsilon}$ with   $f(\partial \Delta) \subset \partial D$, and 
$\mathfrak{Var}_o(f) $ denotes the class of infinitesimal variations of $f$  attached to   $\partial D$ and with  fixed center. \par
\begin{theo}  \label{nicety}
Assume that $D \subset M$ is of the form $D = \{\ \rho < 0\ \}$ for some  $J$-plurisubharmonic  $\rho$  (see \S \ref{plurisubharmonic}, for definition) and that  $\mathfrak{Var}_o(f)$ contains a $(2n-2)$-dimensional $J$-invariant vector space, generated by  infinitesimal variations $e_i$, $J e_i$, $1 \leq i \leq n-1$,   such that   the maps 
 $\z^{-1} \cdot e_i(\z) \ ,\ \z^{-1} \cdot J e_i(\z): \overline \D \longrightarrow TM$ 
 are of class $\cC^{\a, \epsilon}$ on $\overline \D$. Assume also that, for any $\z \in \overline \D$,  the set $\{e_i(\z), J e_i(\z)\} \subset  T_{f(\z)}  M$ span a subspace,   which is complementary to $T_{f(\z)} f(\D)\subset  T_{f(\z)}  M$.   
Then $f$ is critical if and only if it is stationary. 
\end{theo}
\begin{pf} Let $e_0: \overline \D \longrightarrow TM$ be the map defined by $ e_0(\z) = f_*\left(\left.\frac{\partial}{\partial \Re \z}\right|_\z\right)$.
By hypotheses,   the collection 
\beq \label{basis} \left(e_0(\z), J e_0(\z), \z^{-1} \cdot e_1(\z),\dots, \z^{-1}\cdot e_{n-1}(\z), \z^{-1} \cdot  J e_{n-1}(\z)\right)\eeq
is a basis for $T_{f(\z)} M$ for all $\z \in \D$ and  we may consider a system of coordinates $\xi: (x^0, x^1, \dots, x^{2n-2}, x^{2 n-1}) = (z^0, \dots, z^{n-1}) : \cU \longrightarrow \bR^{2n} = \bC^n$ on a neighborhood $\cU$ of $f(\overline \D)$ such that 
\beq \label{adaptedcoordinates}\left.\frac{\partial}{\partial x^{0}}\right|_{f(\z)} = e_0(\z)\ ,\ \  \left.\frac{\partial}{\partial x^{1}}\right|_{f(\z)}  = J (  e_0(\z))\ ,$$
$$\left.\frac{\partial}{\partial x^{2i}}\right|_{f(\z)} = \z^{-1} \cdot e_i(\z)\ ,\ \  \left.\frac{\partial}{\partial x^{2i +1}}\right|_{f(\z)}  = J ( \z^{-1} \cdot e_i(\z))\ \ \text{for all}\ \ \z \in\D \setminus \{0\}\ .\eeq
If we identify $\cU$ with $\xi(\cU) \subset \bC^n$,  we have that $J|_y = \Jst|_y$ for all $y \in f(\D)$ and   the maps $f$,  $e_i$ and $J e_i$ are of the form  
(here, any vector valued map is  denoted  by a pair,  formed by the base point  and the vector  components): 
\beq \label{basis1}f(\z) = (\Re \z, \Im\z, 0, \dots, 0)\ ,$$
\par
$$ e_i(\z) = ((\Re \z, \Im\z, 0, \dots, 0);  (0, \dots, \underset{\text{\tiny $2i$-th place}}{\Re(\z)}, \dots, 0))\ ,$$
\par
$$ J e_i(\z) = ((\Re \z, \Im\z, 0, \dots, 0);  (0, \dots, \underset{\text{\tiny $(2i+1)$-th place}}{\Im(\z)}, \dots, 0))\ .\eeq
A map $W = (f, v^j \frac{\partial}{\partial x^j} ) : \D \longrightarrow T \bR^{2n} \simeq T M$ is an infinitesimal variation (i.e. $\bJ$-holomorphic)  if and only if  it is  solution  of the  p.d.e. system 
\beq \label{pde} \frac{\partial v^i}{\partial \Re \z} + \left.J^i_j\right|_f \frac{\partial v^j}{\partial \Im \z} + 
\left.\frac{\partial J^i_j}{\partial x^k}\right|_f v^k\frac{\partial f^j}{\partial \Im \z} \! = \!  \frac{\partial v^i}{\partial \Re \z} + \Jst{}^i_j \frac{\partial v^{j}}{\partial \Im \z} + 
\left.\frac{\partial J^i_1}{\partial x^k}\right|_f  v^k\! =\! 0. \eeq
By hypotheses,  the fields in  \eqref{basis1} are solutions of  \eqref{pde}. From this and the fact  that   the components $J^i_j|_{f(\D)} =  \Jst{}^i_j|_{f(\D)} $ are  constant on $f(\D)$  we get   that
\beq\label{3.9} \frac{\partial J^i_1}{\partial x^k}=  \frac{\partial J^i_j}{\partial x^0} =  \frac{\partial J^i_j}{\partial x^1} = 0\qquad \text{at all points of} \ f(\D)\ .\eeq
From this and   the explicit formulae in coordinates for $\bJ$ and $\wt \bJ$, one can directly check  that  a map $W = (f, w^j \frac{\partial}{\partial z^j} + \overline w^j \frac{\partial}{\partial \overline z^j}): \D \longrightarrow T \bC^{n} \simeq TM$  (resp.   $\wt f = (f, g_i d z^i + \overline g_i d \overline z^i) : \D \longrightarrow T^* \bC^{n}  \simeq T^*M $) is  $\bJ$-  (resp.  $ \wt \bJ$-) holomorphic  if and only if the functions  $w^j$  (resp. $g_i$) are  holomorphic  in the classical sense. 
Assume that $f$ is stationary and that
$\wt f =   ((\z, 0, \dots, 0); g_i dz^i + \overline g_i d \overline z^i)$  is a stationary lift of $f$ and  $W : \overline \D \longrightarrow TM$ is an infinitesimal variation,  attached to $\partial D$ and with  fixed central direction. Then, for any $\z \neq 0$,  $W(\z)$ is  the form $W(\z) =\z \cdot   \mu^j(\z) \cdot e_i(\z) + \overline{\z \cdot   \mu^j(\z) \cdot e_i (\z)}$, 
 for some holomorphic   $\mu^j: \D\longrightarrow \bC$,  and 
$W_*\left(\left.\frac{\partial}{\partial \Re \z}\right|_{0}\right) =  \mu^j(0) \cdot e_j(0)$.
By Definition \ref{variationalfields} (c), we get that  
$\mu^0(0) \in \bR$ and $\mu^i(0) = 0$ for $i \neq 0$. On the other hand, by Definition \ref{variationalfields}(b), 
  the function 
$$\varphi: \overline \D \longrightarrow \bR\ ,\qquad \varphi(\z) = \wt f(\zeta)\left(\z^{-1}\cdot W(\z)\right) = 
\Re(\mu^j g_j)(\z)\ ,$$
is   so that  $\varphi|_{\partial \D} \equiv 0$. Since $\varphi$ is harmonic,   $\varphi \equiv 0$ and  in particular $\mu^0(0) \Re(g_0(0)) = 0$. By \cite{GJ},  Cor. 2.5,  $\wt f\left(\left.\frac{\partial}{\partial \Re \z}\right|_{0}\right) = g_0(0)  \neq 0$ and  $\mu^0(0) = 0$, i.e.  $f$ is critical. 
Conversely, assume $f$  critical and let  $\wt f: \overline \D\longrightarrow T^* M$ be  defined by
$\wt f(\z) =  (\z, 0\dots, 0; dz^0 + d\overline z^0)$. 
By previous observations, $\wt f$ is $\wt \bJ$-holomorphic and it is an embedding. From the fact that  $J e_0|_{\z}$, $e_i|_\z$, $J e_i|_\z$  span the 
tangent spaces $T_{f(\z)} \partial D$, $\z \in\partial  \D$,  we get   that $\wt f$ is  a stationary lift of $f$ and  that  $f$ is stationary.\end{pf}
\par
\medskip
\subsection{Critical foliations of  normal forms}
Let $(B^n, J)$ be an almost complex circular domain in normal form. We  recall that, for any $v \in S^{2n -1} \subset T_0 B^n$, the straight disk  $f(\z) = \z \cdot v$ is stationary for $(B^n, J)$. \par
For any $w \in T_v S^{2n - 1}$ and any smooth curve $\g^{(w)}_t$  in $S^{2n-1}$ with $\g^{(w)}_0 = v$ and $\dot \g^{(w)}_0 = w$, we may consider the infinitesimal variation, with fixed center and attached to $\partial B^n$,  defined by 
\beq\label{coorvar} W^{(w)}: \overline \D \to T \bC^{2n}\ ,\qquad W^{(w)}(\z) =   \left.\frac{ d\left(\z\cdot \g^{(w)}_t\right)}{dt}\right|_{t =0} =\z \cdot w \in T_{f(\z)} \bC^n\ .\eeq
These maps form a $(2n-2)$-subspace  $\wt{ \mathfrak{Var}}_o(f) $ of the vector space of infinitesimal variations in $\mathfrak{Var}_o(f)$ such that: 
\begin{itemize}
\item[a)]Ê  
 $ \left\{ W(\z) , \ W \in \wt{\mathfrak{Var}}_o(f) \ \right\} = \cH_{f(\z)}$ for any  $\z \neq 0$; 
 \item[b)]  for any $W \in  \wt{\mathfrak{Var}}_o(f)$, the map $\a(\z) \= |\z|^{-1}|W(\z)|$ is constant  on $ \overline\D \setminus \{0\}$.
 \end{itemize}
 Moreover, 
  \begin{lem} \label{lemmetto}Ê The space $J \wt{\mathfrak{Var}}_o(f)$  is included in  $\mathfrak{Var}_o(f)$ if and only if 
 $\cL_{Z^{01}} J = 0$,   with  $Z^{01} = \overline z^i\frac{\partial}{\partial \overline z^i}$.
 \end{lem}
 \begin{pf} Consider an open subset $\cV \subset S^{2n-1} \subset T_0 B^n$ and a field of   real frames  $(e^o_1(v), J_0 e^o_1(v)\dots, J_0 e^o_{n-1}(v))$, $v \in \cV$, for the holomorphic subspaces $H_v \subset T_v S^{2n-1}$ of $S^{2n-1}$. We denote by $e_i(v, \cdot)$,  $(J_0e_i)(v, \cdot)$ the corresponding infinitesimal variations along $f^{(v)}(\z) = \z \cdot v$, i.e.  
 $$e_1(v; \z) \=\z   \cdot  e^o_1(v)\ , \qquad \dots\qquad  ,\quad \ (J_0e_{n-1})(v; \z) \= \z  \cdot  (J_0 e^o_{n-1}(v))\ . $$
Notice that  the points  $ f^{(v)}(\z)$,   with $ (v,\z) \in \cV \times \{ \D \setminus \{0\}\}$,   fill an open subset  $\cU \subset B^n \setminus \{0\}$,  that  the ordered set of vector fields 
 \small{$$\left(e_0( v;\z) \= f^{(v)}_*\left(\left.\frac{\partial}{\partial x}\right|_\z\right), J e_0( v;\z) \= f^{(v)}_*\left(\left.\frac{\partial}{\partial y}\right|_\z\right), e_1(v;\z) ,  \dots, (J_0e_{n-1})(v; \z)\right)
$$}
\par
\noindent
\normalsize is a frame field on $\cU \simeq \cV \times \{ \D \setminus \{0\}\}$ and that  the field   
 $  \frac{1}{2}\left(e_0(v; \z) + iJ e_0(v; \z)\right)$ $=$  $f^{(v)}_*\left(\left.\frac{\partial}{\partial \overline \z}  \right|_\z\right) $ 
is   a generator for  $\cZ^{01}$.  One can also check that
 \beq \label{Zinvariance} \cL_{Z^{01}} e_i = \cL_{Z^{01}} (J_0e_i) = 0\qquad \text{for}\ 1 \leq i \leq n-1\ .   \eeq
 We claim that 
 $\cL_{Z^{01}} J e_i = \cL_{Z^{01}} J (J_0e_i) = 0$ (or, equivalently, that $\cL_{Z^{01}} J|_{\cU} = 0$) if and only if the fields  $J e_i(v; \z)$ and $J (J_0e_i)(v; \z) $ are   in $\mathfrak{Var}_o(f^{(v)})$;   
 by arbitrariness of  $\cV \subset S^{2n-1}$ this conclude the proof.
To check the claim, let us fix a straight disk $f^{(v)}(\z) = \z \cdot v$. By construction, the  fields $\z^{-1} \cdot e_i(v; \z)$ and  $\z^{-1} \cdot (J_0e_i)(v; \z)$  are of class $\cC^{\infty}$ at any  $\z \in \overline \D$. We may therefore consider a system of coordinates $\xi = (z^0 = x^0 + i x^1,  \dots, z^{n-1} = x^{2n - 2}Ê+ i x^{2n-1})$ around  $f^{(v)}(\overline \D)$ satisfying \eqref{adaptedcoordinates}. 
As in the proof of  Theorem \ref{nicety},  we have  that a vector field $W: \D \longrightarrow TB^n $ along $f^{(v)}(\D)$ is $\bJ$-holomorphic if and only if the complex  functions $w^j: \D \longrightarrow  \bC$ such that 
 $$W_{\z} = \Re(w^j(\z))\  e_j(f^{(v)}(\z)) + \Im(w^j(\z)) (J_0e_j)(f^{(v)}(\z))$$
 are holomorphic. By \eqref{Zinvariance}, such holomorphicity condition 
  is equivalent to  $\cL_{Z^{01}}ÊW = 0$. Since, by construction, the fields  $J e_i $ and $J (J_0e_i) $  satisfy (a) and (b) of Definition \ref{variationalfields},  they are  in $\mathfrak{Var}_o(f^{(v)})$ if and only if they are $\bJ$-holomorphic, i.e. if and only if $\cL_{Z^{01}} J e_i = \cL_{Z^{01}} J (J_0e_i) = 0$.
  \end{pf}
Let us introduce the following definition. \par
\begin{definition} Let $(B^n, J)$ be an almost complex domain  of circular type in normal form. We call it  {\it nice} if the distribution  $\cH$ defined in \eqref{2.12} is $J$-invariant.  We call it  {\it very nice} if for any straight disk $f(\z) = \z \cdot v$,  the associated vector space $\wt{\mathfrak{Var}}_o(f)$ is     $J$-invariant. \par
An almost complex domain $(D, J)$ of circular type with center $x_o$  is called {\it nice} (resp. {\it very nice}) if it has a nice (resp. very nice)  normal form. 
\end{definition} 
Motivation for considering such notions comes from  the following\par
\begin{prop} \label{cor} The stationary disks 
of the circular type foliation $\cF^{(x_o)}$ of a  very nice almost complex domain $(D, J)$ of the form $D = \{\ \rho < 0\ \}$ for some  $J$-plurisubharmonic  $\rho$, are critical.
\end{prop} 
\begin{pf} 
With no loss of generality, assume that $(D, J) = (B^n, J)$ is in normal form and that its  foliation of circular type  is given by the straight disks $f(\z) = \z  \cdot v$,  $v \in S^{2n-1} \subset T_0 B^n \simeq \bC^n$. Fix $v \in S^{2n-1}$ and let $H_v \subset T_v S^{2n-1}$ be the holomorphic tangent space   at $v$ and  $(e^o_1, \dots, e^o_{n-1})$   a basis over $\bC$    for $H_v$.  Consider  the infinitesimal variations in $\wt{\mathfrak{Var}}_o(f)$ defined in \eqref{Zinvariance} 
$$e_1(\z) \=\z   \cdot  e^o_1\ , \qquad \dots\qquad  , \ e_{n-1}(\z) \= \z  \cdot e^o_{n-1}\ ,\qquad \z \in \overline \D\ .$$
By construction, the fields $\z^{-1} \cdot e_i$ and  $\z^{-1} \cdot J e_i$, defined at the points of $\overline \D$, are of class $\cC^{\infty}$. Being $(B^n, J)$ very nice, they span  $\cH_z \subset T_{f(\z)} \bC^n$, which is  complementary to $T_{\f(\z)}Êf(\overline \D)$.  By Theorem \ref{nicety}, the conclusion follows. \end{pf}
By definitions  any very nice domain is   nice. The converse is not true, as   next proposition and example show. \par
 \medskip
 \begin{prop} \label{chara1}
Let $(B^n, J)$ be an almost complex domain of circular type in normal form, with $J$ given  by  a deformation tensor $\phi =  \phi^\cH + \phi^{\cH, \cZ}$. 
It is  nice if and only if $\phi^{\cH, \cZ}\equiv 0$, while it is very nice if and only if 
\beq\label{verynicety} \phi^{\cH, \cZ}\equiv 0\qquad \text{and}\qquad  \cL_{Z^{01}} \phi^{\cH} = 0\ .\eeq
\end{prop}
\begin{pf}  It follows  from  definitions and Lemma \ref{lemmetto}.  
\end{pf}
\begin{example} Let  $\phi = \phi^\cH$ be a deformation tensor  in $\Hom( \cH^{01}_z,  \cH^{10}_z)$ at any $z \in B^n$,  which is  non zero only on a relatively compact subset, whose closure does   not contain the origin. By Remark \ref{remark} and Proposition \ref{chara1},   $\phi = \phi^{\cH}$ determines an almost complex structure $J$ such that  $(B^n, J)$ is nice.  On the other hand, by assumptions,   there are  straight disks $f$ with   $\phi^{\cH}|_{f(\D)} \not \equiv 0$ and  with $\phi^{\cH}|_{\cV} \equiv 0$ on some   open subset of $\cV \subset f(\D)$. Due to this,   the equality $ \cL_{Z^{01}} \phi^{\cH} = 0$ is not satisfied and $(B^n, J)$ is not very nice. 
\end{example}
\par
\bigskip
\section{Almost complex Monge-Amp\`ere operators}
\setcounter{equation}{0}
\subsection{Plurisubharmonic functions  and  pseudoconvex manifolds}\label{plurisubharmonic}
  Let $(M, J)$ be an almost complex manifold and  $\O^k(M)$, $k \geq 0$, 
  the space of  $k$-forms of $M$. We denote by  $d^c: \O^k(M) \longrightarrow \O^{k+1}(M)$  the classical $d^c$-operator
  $$d^c \a = (-1)^k (J^* \circ d \circ J^*)(\a)\ ,$$ 
where $J^*$   denotes the usual action  of $J$ on $k$-forms, i.e.  
$J^* \b(v_1, \dots, v_k) \= (-1)^k \b(Jv_1, \dots, J v_k)$
(see e.g. \cite{Bs}). If
%
$J$ is integrable,  it is well known that 
$$d^c = i (\overline \partial - \partial)\ ,\quad \partial \overline \partial = \frac{1}{2i}Êd d^c\ ,\quad d d^c = - d^c d$$
and that $d d^c u$ is a $J$-Hermitian 2-form   for any  $\cC^2$-function $u$. Unfortunately, 
{\it when  $J$ is not integrable,  $d^c d \neq -d d^c$  and the 2-forms $d d^c u$, with  $u \in \cC^2(M)$,  are usually  not   $J$-Hermitian}.   In  fact, one has that   
\beq \label{Jinv}  d d^c u(J X_1, X_2) +  d d^c u(X_1, J X_2)  =  4 N_{X_1 X_2} (u)\ ,\eeq
where $N_{X_1 X_2}$ is  the Nijenhuis tensor evaluated  on $X_1$, $X_2$ and is in general non zero. This fact suggests   the following definition.
 \par
\begin{definition} Let  $u: \cU \subset M \longrightarrow \bR$ be of class   $\cC^2$. We call   {\it $J$-Hessian  of $u$ at $x$\/} the symmetric form  $\HH(u)_{x} \in S^2 T_{x} M$, whose associated   quadratic form is $ \cL(u)_{ x}(v)= d d^c  u(v, J v)_{x}$.
By polarization formula and \eqref{Jinv},   one has that, for any $v$, $w \in T_x M$, 
\beq \HH(u)_{x} (v, w) = \left.\frac{1}{2} \left(d d^c u(v,J w ) + dd^c u(w,J v)\right)\right|_{x} = $$
$$ =   d d^c u(v,J w )_{x} - 2 N_{v w}(u) \ .  \eeq
\end{definition}
\smallskip
We remark  that $\HH(u)_x$ is not only symmetric, but  also $J$-Hermitian, i.e. 
$\HH(u)_{x}(Jv, Jw) = \HH(u)_{x}(v, w) $ for any $v, w$. It is therefore  associated with  the Hermitian antisymmetric tensor
\beq \HH(u)(J\cdot , \cdot ) = \frac{1}{2} \left(dd^c u(\cdot ,\cdot) + d d^c u(J\cdot ,J \cdot )\right) =  \frac{1}{2} \left(dd^c u + J^* d d^c u\right)  .\eeq
The quadratic form $ \cL(u)_{ x}(v) =   d d^c u(v, Jv)|_{x}$
is the so-called {\it Levi form of $u$ at $x$\/} (see e.g. \cite{CGS2}) and it is  tightly  related with the notion  of $J$-plurisubharmonicity.   On this regard, we recall that  an upper semicontinuous function $u: \cU \subset M \longrightarrow \bR$  is called {\it $J$-plurisubharmonic\/}Ê   if,   for any $J$-holomorphic disk $f : \D \longrightarrow \cU \subset M$,  the composition $u \circ f: \D \longrightarrow \bR$ is subharmonic. 
By   simple arguments (similar to those used for  complex manifolds),  whenever  $u$ is in $\cC^2(\cU)$   one has that  {\it $u$ is $J$-plurisubharmonic   if and only if $ \cL(u)_{ x}(v) = \HH(u)_{x}(v,v) \geq  0$ for any $x \in \cU$ and  $v \in T_xM$\/}.\par
\medskip
This motivates the following  generalizations of  classical notions (see e.g. \cite{De1}).
 In the following, for any  $\cU \subset M$,  the symbol  $\Psh(\cU)$ denotes the class  of  $J$-plurisubharmonic functions on $\cU$.\par
\begin{definition} Let $(M, J)$ be an almost complex manifold and $\cU \subset M$ an open subset.  We say that $u \in \Psh(\cU)$ is {\it strictly $J$-plurisubharmonic}   if: 
\begin{itemize}
\item[a)]  $u \in L^1_{\operatorname{loc}}(\cU)$; 
\item[b)]   for any $x_o \in \cU$ there exists a neighborhood $\cV$ of $x_o$ and $v \in \cC^2(\cV) \cap \Psh(\cV)$ for which 
$\HH(v)_{x}$ is positive definite at all points  and $u - v$ is in $\Psh(\cV)$. 
\end{itemize}
In particular, $u \in \Psh(\cU) \cap \cC^2(\cU)$ is strictly plurisubharmonic if and only if $\HH(u)_x$ is positive definite at any $x \in \cU$.\par
 The almost complex manifold $(M, J)$  is called  {\it weakly} (resp. {\it strongly}) {\it pseudoconvex} if it admits  a $\cC^2$  exhaustion $\t: M \longrightarrow ]- \infty, \infty[ $, which is plurisubharmonic (resp. strictly plurisubharmonic) (\footnote{Strongly pseudoconvex manifolds are called {\it almost complex Stein manifolds} in  \cite{DS}.}).
\end{definition}
\medskip
\subsection{Maximal plurisubharmonic functions}
 $J$-plurisubharmonic functions share most of the basic properties of classical plurisubharmonic functions. For instance, for any open domain  $\cU \subset M$,  the  class  $\Psh(\cU)$ is a convex cone and  a lattice, as  for  domains in complex manifolds. In fact, given $u_i \in \Psh(\cU)$ and  $\l_i \in \bR$, also  the functions
$u = \sum_{i = 1}^n \l_i u_i$ and $u' = \max \{\ u_1, \dots, u_n\ \} $
are in $\Psh(\cU)$. \par
It is therefore natural  to  consider   the following  notion of   ``maximal''  $J$-plurisubharmonic functions. This  and  next theorem  indicate which operator  should be  considered as natural generalization of  the classical  complex Monge-Amp\`ere operator.  \par
\begin{definition}  Let $D$ be a domain in a strongly pseudoconvex almost complex manifold $(M, J)$. 
  A function $u \in  \Psh(D)$ is called {\it maximal} if for any open $\cU \subset \subset D$
  and  $h \in \Psh(\cU)$ satisfying the condition 
\beq \label{condmax1}  \limsup_{z \to x} h(z) \leq u(x)\quad \ \ \text{for all} \ x \in \partial \cU\ ,\eeq
  one has that $h \leq u|_{\cU}$.   
\end{definition}
\begin{theo} \label{firsttheorem} Let $D \subset M$ be a domain of a strongly pseudoconvex almost complex manifold $(M, J)$ of dimension $2n$. A function 
$u \in \Psh(D) \cap \cC^2(D)$  is maximal if and only if it satisfies
\beq \label{MA1} \left(d d^c u + J^*(d d^c u)\right)^n = 0\ . \eeq
\end{theo}
\begin{pf} 
Let  $\t: M \longrightarrow ]-\infty, + \infty[$ be  a $\cC^2$ strictly plurisubharmonic exhaustion  for $M$ and assume that $u$ satisfies \eqref{MA1}. We  need to  show that for any $h \in \Psh(\cU)$ on an  $\cU \subset \subset D$  that 
satisfies \eqref{condmax1}, one has that  $h \leq  u|_{\cU}$. Suppose not and pick $\cU \subset \subset D$   and $h \in \Psh(\cU)$, so that \eqref{condmax1} is true but  there exists 
$x_o \in \cU$ with $ u(x_o) < h(x_o)$.  Let 
$\l > 0$ so small that 
$$h(x_o) + \l\left( \t(x_o)- M \right) > u(x_o)\ ,\qquad \text{where} \ M = \max_{y\in \overline{\cU}} \t(y)\ ,$$
and  denote by $\wh h$ the function 
\beq \wh h \=\left. h + \l (\t -  M )\right|_{\cU}\ . \eeq
By construction,  $\wh h \in \Psh(\cU) $,  satisfies \eqref{condmax1} and $(\wh h - u)(x_o)  > 0$. In particular, $\wh h - u$ achieves its maximum at some  inner point $y_o \in \cU$. 
 Now, we remark that  \eqref{MA1} is equivalent to say  that, for  any $x \in D$,  there exists $0 \neq v \in T_x M$ so that 
\beq \label{MA2} \left(d d^c u + J^*(d d^c u)\right)_x(v , Jv) = \HH_x(u)(v, v) = 0\ .\eeq
Let  $0 \neq v_o \in T_{y_o} M$ be a vector for which \eqref{MA2} is true  and let  $f: \D \longrightarrow M$ be  a $J$-holomorphic disk so that  $f(0) = y_o$ and  with 
$$f_*\left(\left.\frac{\partial}{\partial x}\right|_{0}\right) = v_o \ ,\qquad f_*\left(\left.\frac{\partial}{\partial y}\right|_0\right) = 
f_*\left(\Jst \left.\frac{\partial}{\partial x}\right|_{0}\right) = Jv_o\ .$$
Then, consider the  function $G : \D \longrightarrow \bR$ defined by 
\beq G \= \wh h \circ f - u \circ f =  h\circ f + (\l \t - \l M - u) \circ f \ .\eeq
We claim that there exists a disk $\D_r =\{|\zeta| < r \}$ such that $G|_{\D_r}$ is subharmonic. 
In fact, since $\t$ is $\cC^2$ and strictly plurisubharmonic and  $\HH(u)_{y_o}(v_o,v_o) = 0$, we have that 
$$0 < \HH((\l \t - \l M - u))_{y_o}(v_o,v_o)  =  2i \left.\partial \overline \partial ((\l \t - \l M - u) \circ f )\right|_0\ .  $$
Hence, by continuity, there exists $r > 0$ so that 
$$0 <  2 i \left.\partial \overline \partial ((\l \t - \l M - u) \circ f )\right|_\z\qquad \text{for any } \ \ \z \in \overline{\D_r}\ .  $$
It follows that  $(\l \t - \l M - u) \circ f |_{\D_r}$ is strictly subharmonic  and that  $G|_{\D_r}$ is subharmonic,  being  sum of  subharmonic functions.   At this point, it suffices to observe that, since $y_o$ is a point of maximum for $\wh h - u$ on $f(\D) \subset \cU$, 
then $0 = f^{-1}(y_o) \in \D_r$ is an inner point of maximum for $G|_{\D_r}$. In fact, from this and  the maximum principle,  we get  that $G|_{\D_r}$ is constant and hence that  $h\circ f|_{\D_r}$ is $\cC^2$ 
with $2 i  \left.\partial \overline \partial(h\circ f) \right|_{\D_r}< 0$,   contradicting the hypothesis  on subharmonicity of  $h\circ f$. \par
\smallskip
Conversely, assume that $u \in \cC^2(D) \cap \Psh(D)$ is maximal, but  that \eqref{MA1} is not satisfied, i.e. that  there exists $y_o\in D$  for which   $\HH_{y_o}(u)(v, v) > 0$ for any $0 \neq v \in T_{y_o}ÊM$. By  Lemma  \ref{small}, there exist a relatively compact neighborhood $\cU$ of $y_o$  and a $(J, J')$-biholomorphism between $(\cU, J)$ and  $(B^n, J')$, with $J'$ arbitrarily close in $\cC^2$ norm to the standard complex structure. Due to this,   we may assume that  $\t =$ $\t_o \circ \varphi$, with $\t_o(z) = |z|^2$,  is a $\cC^2$ strictly $J$-pluri\-sub\-harmonic  ex\-haustion on $\cU$, tending to  $1$ at the points of $\partial \cU$.  Hence, there is   a constant $c > 0$ such that
$$\HH_{x}(u + c(1 -   \t))(v, v) = \HH_x(u)(v, v) - c \HH_{x}(\t)(v,v) \geq 0\ ,$$
for all $x \in \cU$ and $v \in  T_x M \simeq \bR^{2n}$  with $|v| = 1$.
This means that 
$$\wh h \= \left.\left(u + c(1 -  \t)\right)\right|_{B_{y_o}(r)} $$
 is in $\cC^2(\cU) \cap \Psh(\cU)$,  satisfies \eqref{condmax1} and,  by maximality of $u$,  satisfies $\wh h \leq u $ at all points of $\cU$. 
But there is also an  $\epsilon > 0$  such that  $ \emptyset \neq  \t^{-1}([0, 1-\epsilon[) \subsetneq \cU$ and hence such that, on  this subset,   $\wh h \geq  u + c \epsilon  > u$,  contradicting the maximality of $u$.
\end{pf}
\subsection{Green functions of nice circular domains} The results of previous section show that \eqref{MA1} is a natural   analogue of classical complex Monge-Amp\`ere equation for  domains in $\bC^n$ and that the solutions of \eqref{MA1} are interesting  biholomorphic invariants of strongly pseudoconvex domains. This motivates the following generalized notion of Green functions (see e.g. \cite{Be}).\par
\begin{definition} \label{Green}ÊLet $D$ be a domain in a strongly pseudoconvex,  almost complex manifold $(M, J)$. We call {\it almost pluricomplex Green function with pole at $x_o \in D$} an exhaustion $u: \overline{D} \longrightarrow [- \infty, 0]$  such that
\begin{itemize}
\item[i)] $u|_{\partial D} = 0$ and  $u(x) \simeq \log \|x - x_o\|$ when $x \to x_o$,  for some Euclidean metric $\|\cdot \|$ on a neighborhood of $x_o$; 
\item[ii)] it is $J$-plurisubharmonic; 
\item[iii)] it is a solution of the generalized Monge-Ampere equation $ \left(d d^c u + J^*(d d^c u)\right)^n = 0$ on $D \setminus\{x_o\}$.
\end{itemize}
\end{definition}
Notice that,  if a Green function with pole $x_o$ exists,  by a direct consequence of property of maximality (Theorem \ref{firsttheorem}) it is unique. \par
\smallskip 
Consider now an almost complex  domain $D$ of circular type in $(M, J)$  with center $x_o$. Denoting by $\exp: \wt B^n \longrightarrow \wt D$ the corresponding Riemann map, we call
 {\it standard exhaustion of $D$} the map
$$\tau_{(x_o)}: D \longrightarrow [0, 1[\ ,\qquad \tau(x) = \left\{\begin{matrix} |\exp^{-1}(x)|^2 & \text{if}\ x \neq 0\ ,\\
\phantom{a} & \ \\
0 & \text{if}\ x = x_o\ .\end{matrix} \right.$$
When  $D$ is  in normal form, i.e. $ D = (B^n, J)$ with $J$ almost $L$-complex structure, its standard exhaustion is just  $\tau_o(z) = |z|^2$.
\par
\begin{prop} \label{lastprop}ÊLet $D$ be a  domain of circular type in  $(M, J)$ with center $x_o$ and standard exhaustion $\t_{(x_o)}$. If    $u = \log \tau_{(x_o)}$ is $J$-plurisubharmonic, then $u$ is  an almost pluricomplex Green function with pole at $x_o$. 
\end{prop}
\begin{pf} With no loss of generality, we may assume that  the domain  is in normal form, i.e. $ D = (B^n, J)$ and  $\t_{(x_o)}(z) = \t_o(x) = |x|^2$. Since $\t_o$ is smooth on $B^n \setminus \{0\}$ and  $u = \log \t_o$ is  $J$-plurisubharmonic,  we have  that  $\HH(u)_x \geq 0$ for any $x \neq 0$. On the other hand, for any straight disk 
$f: \D \longrightarrow B^n$ of the form  $ f(\z) = v \cdot \z$, we have that 
 $u \circ f$ is harmonic and $\HH(u)_{f(\z)}(v,v) = 0$ for any $\z \neq 0$. This means that    $\HH(u)_x \geq 0$ has at least one vanishing eigenvalue at any  point of $B^n \setminus \{0\}$ and means that  \eqref{MA1} is satisfied. Other conditions of Definition \ref{Green}   can be checked directly. 
\end{pf}
When $J$ is integrable, {\it the standard exhaustion $u = \log \tau_{(x_o)}$ of  the normal form of a  
domain of circular type is automatically  plurisubharmonic }  (\cite{PS}). In the almost complex case, this is no longer true, as the following example shows. 
\par
\begin{example}\label{example58}  
Consider a quadruple of vector fields  $(Z, \Jst Z, E, \Jst E)$  on $\wt B^2$, 
determined as follows. The field $Z$ has been defined in \S \ref{domainsofcirculartype} and, in coordinates \eqref{interestingcoord},   is of the form  $Z_z = \Re\left( z^0\left.\frac{\partial}{\partial z^0}\right|_z\right)$ at any  $z \in \wt B^2 \setminus \pi^{-1}(0)$. 
The  field $E$   is any vector field in the distribution $\cH$ that satisfies the conditions
\beq [Z, E] = [\Jst Z, E] = 0\ ,\qquad [E, \Jst E] = - \Jst Z \ .\eeq
It is uniquely determined, up to a smooth family of  unitary transformations of the subspaces $\cH_z \subset T_z \wt B^2$.  
Notice that the standard holomorphic bundle $T^{10} \wt B^2$  is  generated at all points by the complex vector fields $Z^{10} = Z - i \Jst Z$,  $E^{10} = E - i \Jst E$. In the following, we denote by $(E^{10*}, E^{01*}, Z^{10*}, Z^{01*})$ the field of complex coframes,  which is dual to $(E^{10}, E^{01} = \overline{E^{10}}, Z^{10}, Z^{01}= \overline{Z^{10}})$ at all points.\par
Consider a deformation tensor $\phi \in \Hom(\cH^{01}, \cZ^{10} + \cH^{10})$ of the form 
$\phi_z =  h(z) Z^{10}_z \otimes E^{01}{}^*_z$ for some 
 smooth real valued function $h: \wt B^n \longrightarrow \bR$, which  is constant on all spheres $S_c = \{\ \tau_o(z) = c\ \}$ (i.e. $h = h(|z^0|)$)  and is equal to $0$ on 
an open neighborhood of $\pi^{-1}(0) = \bC P^1$. \par
 By definitions and Remark \ref{remark}, the deformation tensor $\phi$  determines an  almost complex structure $J$, with  $J$-holomorphic spaces $T^{10}_{Jz} \wt B^n = \bC Z^{10}_z \oplus \bC\wt E^{10}_z$,   $\wt E^{10}_z \=  E^{10}_z + h(z) Z^{01}_z$,  that satisfies (i), (iii) and (iv) of Definition \ref{Lstructure}. We claim that the system \eqref{diffsyst} is alway solvable, so that $J$ is an almost L-complex structure and  $(\wt B^n, J)$ is an almost complex domain of circular type in normal form. This  claim can be checked observing that  the  components of $J$ in coordinates \eqref{interestingcoord}  along the disk $f(\z) = \z \cdot v$ with $v = (0,1)$, are  such that 
 $$J^1_1 = i  = - J^{\bar 1}_{\bar 1}\ ,\ \  J^1_{\bar 1} = J^{\bar 1}_{1} = 0$$
 and $J_A^0 = J_A^0(|z^0|)$, $J_A^{\bar 0} = J_A^{\bar 0}(|z^0|)$ for any $A = 1, \bar 1$. So, \eqref{diffsyst} reduces to 
 \beq \label{diffsystter} \left\{ \begin{array}{lr}  g_{1}{},_{ \bar\z}  = \bF_1
& \text{on}\  \D\ ,\\
 \ & \ \\
g_1 = k (\z^2 - 1) & \!\!\! \text{on}\ \  \partial \D \ ,\end{array}\right.\eeq
 with $\bF_1: \overline \D \longrightarrow \bC$ depending only on $\rho = |\z|$ and $k$ constant. Consider the map
 $$\wt \bF_1: \overline \D \longrightarrow \bC\ ,\qquad \wt \bF_1 \= - \frac{1}{\pi} \int_{\overline \D} \frac{\bF_1(|w|)}{w - \z} d w \wedge d \bar w\ .$$
It is  such that 
 $$\wt \bF_1{},_{\bar \z} =  \bF_1 \quad \text{and}\quad  \int_{\partial \D} \wt \bF_1 \z^n d \z =  2 i 
 \int_{\overline \D} \bF_1(\r) \r^{n+1} e^{i n \vartheta} d\rho\wedge d \vartheta  = 0$$
 for any $n \geq 0$. Hence, if  $h_1: \overline \D \longrightarrow \bC$ is a  holomorphic map such that $h_1|_{\partial \D} =    \wt \bF_1|_{\partial \D}$,  the map $g_1 = \wt \bF_1 - h_1 +   k(\z^2 - 1)$ is a  solution to \eqref{diffsystter}. Since the solvability of \eqref{diffsyst} is   independent on the choice of coordinates, this concludes the proof of the claim.
 \par
\medskip
Now, we want to show that  {\it if  $h \not \equiv 0$, the function  $u = \log \t_o$  is not $J$-plurisubharmonic}. For this,   we first observe  that,  for any  pair of real  vector fields $X, Y$, if we set $X^{10} = X - i J X$, $Y^{01}Ê= Y + i JY$
\beq \label{hessian1} \HH(u)(X, Y) =  \frac{1}{2} \Im d d^c u(X^{10}, Y^{01}) = $$
$$ =   \frac{1}{2}\Im \left( i X^{10}(Y^{01}(u)) + i Y^{01}(X^{10}(u)) + J [X^{10}, Y^{01}] (u)\right)\eeq
We also recall that  (here, 
$(\cdot)_Z \= Z(\cdot)$ denotes derivation along $Z$)
\beq Z^{10}(u) = 1 = Z^{01}(u)\ , \quad E^{10}(u) = 0\ ,  \quad \wt E^{10} (u) = h = \wt E^{01}(u)\ ,\eeq
\beq \wt E^{10}(\wt E^{01}(u)) =  h h_Z = \wt E^{01}(\wt E^{10}(u))\ , \eeq
\beq [E^{10}, E^{01}] = - 2 i [E, JE]  = -  i 2 J Z\ ,\eeq
\beq  J [\wt E^{10}, \wt E^{01}] (u) =  J [E^{10},  E^{01}] (u) + i h h_Z (Z^{10} + Z^{01})(u) =  i 2 (1 +  h h_Z) \ ,\eeq
\beq \label{hessian2}  J [\wt E^{10}, Z^{01}] (u)  = i h_Z Z^{01}(u) = i h_Z\ .\eeq
Hence, if we set   $\wt E =  \Re(\wt E^{10})$ and $J \wt E =  \Re(i \wt E^{10})$,  using \eqref{hessian1} - \eqref{hessian2} we may conclude that   
$$\HH(\wt E, \wt E) = 1 + 2 h h_Z\ ,\qquad \HH(\wt E, J \wt E) = 0\ ,$$
$$
\HH(\wt E, Z) = h_Z\ ,\quad \HH(\wt E, J Z) = \HH(Z, Z)  = \HH(Z, JZ) = 0\ , $$
so that  the matrix $H$,  with entries given by the components of  $\HH(u)_z$  in the basis   $\cB = (e_1 = \wt E_z, e_2 =  J\wt E_z , e_3 = Z_z, e_4 = JZ_z )$, is 
$$H  =  \left(\begin{matrix}   1 + 2 h h_Z & 0 &  h_Z & 0\\
0 &  1 + 2 h h_Z  & 0 &  h_Z\\
 h_Z & 0 & 0 & 0\\
0 &  h_Z & 0 & 0
\end{matrix}
\right)\ . $$
Since the eigenvalues of $H$  are  $\l_{\pm} = \frac{(1 + 2 h h_Z) \pm \sqrt{(1 + 2 h h_Z)^2 + 4 h_Z^2}}{2}$,  we conclude   that $u$ is $J$-plurisubharmonic if and only if $h_Z \equiv 0$, i.e.   if and only if $h \equiv 0$. 
\end{example}
From previous example, we see that  the standard exhaustion $\t_{(x_o)}$ of an almost complex domain $(D, J)$ of circular type is {\it in general}  not a Green function, independently on how $J$ is  close to an integrable complex structure. However, the property remains valid if one restricts to the class of  nice domains  and to small deformations of integrable structures, as it is shown in next theorem. This property nicely relates to Thm. 6.4  \cite{GJ}  (see also Remark \ref{remark3.2}) on the existence of   extremal disks for domains with small deformations of an integrable complex structure. \par 
\begin{theo} Let $D$ be a nice circular  domain  with standard exhaustion $\t_{(x_o)}$ and normal form $(B^n, J)$. If $J$ is a sufficiently small $\cC^1$-deformation  of  $\Jst$,  then $u = \log \t_{(x_o)}$ is the Green function with pole at $x_o$. 
\end{theo}
\begin{pf} By Proposition \ref{lastprop}, it suffices to   show that, when   $J$ is sufficiently close to $\Jst$, then  $u = \log  |z|^2$    is $J$-pseudoconvex on $B^n \setminus \{0\}$.  For  this, we first claim  that if $(B^n, J)$ is nice, then $\HH(u)_z(\cZ,\cH) = 0$ at any $z \neq 0$. Since $\cZ$ and $\cH$ are both $J$-invariant, this is equivalent to  claim that,  that for any $z \neq 0$, 
\beq  \label{pippo} \left(d d^c u + J^*(d d^c u)\right)_z(Z, X) = 0\ ,\qquad Z  \= \Re\left(z^i \frac{\partial}{\partial z^i}\right)\ ,\ \ X \in \cH\ .\eeq
But this follows from the  fact that $X$ and $J X$ are tangent to the level sets of $u = \log \t_o$,  that $[\cZ, \cH] \subset \cH$ and hence  that
$$d d^c  u( Z, X) =  - Z (J X(u))  + X(J Z(u)) + (J[Z, X]) (u) = X(\Jst Z(u)) = 0$$
$$J^* d d^c  u( Z,  X)\! =\!   J Z ( X(u))  -  J X(Z(u)) + (J[J Z, J X]) (u) = - J X(Z(u)) = \!  0.$$
Since  $\HH(u)|_{\cZ \times \cZ} = 0$,  it remains to check that  $\HH(u)|_{\cH_z \times \cH_z} \geq 0$. Since any sphere $S^{2n-1}_c = \{\ z \in B^n\ , \ \t_o(z) = c \ \}$, $c \in ]0,1]$, is strongly $\Jst$-pseudoconvex  and its (real) holomorphic tangent distribution  distribution is $\cH|_{S^{2n-1}_c}$,  if $J$ is sufficiently close to $\Jst$, the sphere $S^{2n-1}_c$ is strongly pseudoconvex also w.r.t. $J$ by a continuity argument. This is the same of saying  that   $\HH(u)|_{\cH_z \times \cH_z} \geq 0$ for any $z \neq 0$,  as needed. 
\end{pf}
\par
\medskip
\subsection{Concluding remarks}
It is well known  that on (integrable!) complex manifolds there is a tight connection between the existence of regular plurisubharmonic solutions $u$ of maximal rank for homogeneous complex Monge-Amp\`ere equations and existence of foliations by Riemann surfaces.  In fact,  such foliation is given by  the complex curves along which the solution $u$ is harmonic. \par
This idea was exploited by Lempert in his work on strictly convex domains in $\bC^n$ (\cite{Le}). In this case,  the foliation is made   of the extremal disks for the Kobayashi metric through a fixed point $x_o$, which coincide with stationary  disks through $x_o$. The function $u$ is the pull-back on each disk of the standard Green function with pole at $0$ of  $\D \subset \bC$, 
 so that any  strictly convex domain is a domain of circular type. Conversely, the singular foliation by holomorphic disks, determined by the exhaustion $u$ of a domain of circular type,  is a foliation by Kobayashi extremal disks (\cite{Pt3}). 
 \par
Also for almost complex domains, the plurisubharmonic $\cC^2$ functions, which are harmonic along the leaves of a foliation of circular type $\cF^{(x_o)}$, are  solutions of an almost complex Monge-Amp\`ere equation: \par
\begin{prop} Let $D \subset M$ be a domain  in a strongly pseudoconvex manifold  $(M, J)$ with  a foliation of circular type $\cF^{(x_o)}$  of  $(D, x_o)$.   Let also  $u : D \longrightarrow ]-\infty, +\infty[$ be a function which 
 is in $\Psh(D) \cap \cC^2(D \setminus \{x_o\})$ and so that 
 $u \circ f: \D\setminus \{0\} \longrightarrow \bR$  is harmonic  for any $f \in \cF^{(x_o)}$. Then $u$  is a solution of the generalized Monge-Amp\`ere equation \eqref{MA1} su $D \setminus \{x_o\}$.
\end{prop}
\begin{pf} First of all, we claim that $u$ satisfies the generalized Monge-Ampere equation \eqref{MA1} at all points of $D \setminus \{x_o\}$. In fact, by  (ii) of Definition \ref{foliation},  for any $y \in D$ we know that there  exists a  disk $f: \D \longrightarrow D$ in $\cF^{(x_o)}$ 
with $y = f(\z)$ for some $\z \in \D$. Since $u \circ f$ is harmonic, if we denote by $v = f_*\left(\frac{\partial}{\partial (\Re \z)}\right)$, we have that 
\beq H(u)_y(v,v) = (d d^c u + J^*(d d^c u))_y(v, J v) = \D(u \circ f)_{\z} = 0\ , \eeq
 from which it follows immediately  that $ (d d^c u + J^*(d d^c u))^n_y = 0$. The conclusion follows immediately from  Theorem \ref{firsttheorem}.\end{pf}\par
 Notice that, by Example \ref{example58}, the previous  remark is not so useful to determine almost pluricomplex Green functions,  since in general the function  $u = \log \t_{(x_o)}$, which is determined by the geometric construction, is  not plurisubharmonic. 
 Furthermore, we know  that, in general, stationary disks  of an almost complex domain $D$  are  not extremal disks for the Kobayashi metric of $D$ (\cite{GJ}). 
 Nevertheless, the above properties of strictly convex domains of $\bC^n$ remain valid in a large class of almost complex domains, as it is illustrated in the following theorem, which is direct consequence of  our results.\par
\begin{theo} Let $D$ be an almost complex  domain of circular type with center $x_o$  in   $(M, J)$ strongly pseudoconvex.   If  the  normal form $(B^n, J')$ of $(D, J)$   is  very nice with $J'$ sufficiently close to $\Jst$, then 
\begin{itemize}
\item[a)] the stationary foliation $\cF^{(x_o)}$ consists of extremal disks w.r.t. Kobayashi metric; 
\item[b)] the function  $u = \log \t_{(x_o)}$ is the  almost pluricomplex Green function of $D$ with pole $x_o$; 
\item[c)] the distribution $\cZ_z = \ker \HH(u)_z$ is integrable and the closures of  its integral leaves are the disks in $\cF^{(x_o)}$.
\end{itemize}
\end{theo}

\end{document}